\documentclass[12pt]{amsart}
\usepackage{diagrams}
\diagramstyle[scriptlabels,Postscript=dvips]
\newcommand{\nc}{\newcommand}

\nc{\one}{\mbox{\bf 1}}
\nc{\invtensor}{\underset{\leftarrow}{\otimes}}
\nc{\ad}{\operatorname{ad}}
\nc{\st}{\operatorname{st}}
\nc{\tr}{\operatorname{tr}}
\nc{\str}{\operatorname{str}}
\nc{\tp}{\operatorname{top}}
\nc{\rank}{\operatorname{rank}}
\nc{\corank}{\operatorname{corank}}
\nc{\codim}{\operatorname{codim}}

\nc{\Ext}{\operatorname{Ext}}
\nc{\Sym}{\operatorname{Sym}}
\nc{\sym}{\operatorname{sym}}
\nc{\id}{\operatorname{id}}
\nc{\Id}{\operatorname{Id}}

\nc{\Ree}{\operatorname{Re}}

\nc{\htt}{\operatorname{ht}}
\nc{\Ker}{\operatorname{Ker}}
\nc{\rker}{\operatorname{rKer}}
\nc{\im}{\operatorname{Im}}
\nc{\osp}{\mathfrak{osp}}
\nc{\sgn}{\operatorname{sgn}}
\nc{\F}{\operatorname{F}}
\nc{\Mod}{\operatorname{Mod}}
\nc{\Mat}{\operatorname{Mat}}
\nc{\Soc}{\operatorname{Soc}}
\nc{\Inj}{\operatorname{Inj}}
\nc{\Hom}{\operatorname{Hom}}
\nc{\End}{\operatorname{End}}
\nc{\supp}{\operatorname{supp}}
\nc{\Card}{\operatorname{Card}}
\nc{\Ann}{\operatorname{Ann}}
\nc{\Ind}{\operatorname{Ind}}
\nc{\Coind}{\operatorname{Coind}}
\nc{\wt}{\operatorname{wt}}
\nc{\ch}{\operatorname{ch}}
\nc{\Stab}{\operatorname{Stab}}
\nc{\Sch}{{\mathcal S}\mbox{\em ch}}
\nc{\Irr}{\operatorname{Irr}}
\nc{\Spec}{\operatorname{Spec}}
\nc{\Prim}{\operatorname{Prim}}
\nc{\Aut}{\operatorname{Aut}}
\nc{\Fract}{\operatorname{Fract}}
\nc{\gr}{\operatorname{gr}}

\nc{\loc}{\operatorname{loc}}

\nc{\HC}{\operatorname{HC}}
\nc{\Kl}{\operatorname{Kl}}

\nc{\rDelta}{r\Delta}
\nc{\wdchi}{\widetilde{\chi}}
\nc{\wdH}{\widetilde{H}}
\nc{\wdN}{\widetilde{N}}
\nc{\wdM}{\widetilde{M}}
\nc{\wdO}{\widetilde{O}}
\nc{\wdR}{\widetilde{R}} 
\nc{\wdS}{\widetilde{S}}
\nc{\wdV}{\widetilde{V}}

\nc{\wdC}{\widetilde{C}}

\nc{\fg}{\mathfrak{g}}
\nc{\fl}{\mathfrak{l}}

\nc{\fhg}{\hat{\fg}}
\nc{\fhn}{\hat{\fn}}
\nc{\fhh}{\hat{\fh}}
\nc{\fhb}{\hat{\fb}}
\nc{\hrho}{\hat{\rho}}

\nc{\Ob}{\operatorname{\mathcal Ob}}
\nc{\Dglie}{\operatorname{{\mathcal D}glie}}
\nc{\Fin}{\operatorname{{\mathcal F}in}}

\nc{\Sg}{{\cS(\fg)}}
\nc{\Shg}{{\cS(\fhg)}}
\nc{\Ug}{{\cU(\fg)}}
\nc{\Uhg}{{\cU(\fhg)}}
\nc{\Sh}{{\cS(\fh)}}
\nc{\Uh}{{\cU(\fh)}}
\nc{\Uhh}{{\cU(\fhh)}}
\nc{\Zg}{{{\mathcal{Z}}(\fg)}}

\nc{\tZg}{{\widetilde{\mathcal Z}({\mathfrak g})}}
\nc{\Zk}{{\mathcal Z}({\mathfrak k})}

\nc{\Up}{{\mathcal U}({\mathfrak p})}
\nc{\Ah}{{\mathcal A}({\mathfrak h})}
\nc{\Ag}{{\mathcal A}({\mathfrak g})}
\nc{\Ap}{{\mathcal A}({\mathfrak p})}
\nc{\Zp}{{\mathcal Z}({\mathfrak p})}
\nc{\cZ}{\mathcal Z}
\nc{\cS}{\mathcal S}
\nc{\cO}{\mathcal O}
\nc{\cA}{\mathcal A}
\nc{\cU}{\mathcal U}
\nc{\cH}{\mathcal H}
\nc{\cM}{\mathcal M}
\nc{\cL}{\mathcal L}
\nc{\cF}{\mathcal F}

\nc{\fo}{\mathfrak o}

\nc{\CO}{\mathcal O}
\nc{\Cl}{\mathcal {C}\ell}

\nc{\cN}{\mathcal N}
\nc{\cB}{\mathcal B}
\nc{\cY}{\mathcal Y}
\nc{\zq}{\mathpzc q}

\nc{\fn}{\mathfrak n}
\nc{\fm}{\mathfrak m}
\nc{\fp}{\mathfrak p}
\nc{\fs}{\mathfrak s}
\nc{\fh}{\mathfrak h}
\nc{\ft}{\mathfrak t}
\nc{\fk}{\mathfrak k}
\nc{\fb}{\mathfrak b}
\nc{\fI}{\mathfrak I}
\nc{\veps}{\varepsilon}

\nc{\fsl}{\mathfrak{sl}}
\nc{\fgl}{\mathfrak{gl}}
\nc{\fso}{\mathfrak{so}}
\nc{\fpgl}{\mathfrak{pgl}}
\nc{\fpq}{\mathfrak{pq}}
\nc{\fq}{\mathfrak q}
\nc{\fsq}{\mathfrak{sq}}
\nc{\fpsq}{\mathfrak{psq}}

\nc{\fpo}{\mathfrak{po}}
\nc{\dirlim}{\underset{\rightarrow}{\lim}\,}
\nc{\nen}{\newenvironment}
\nc{\ol}{\overline}
\nc{\ul}{\underline}
\nc{\ra}{\rightarrow}
\nc{\lra}{\longrightarrow}
\nc{\Lra}{\Longrightarrow}

\nc{\Lla}{\Longleftarrow}

\nc{\Llra}{\Longleftrightarrow}

\nc{\thla}{\twoheadleftarrow}

\nc{\hra}{\hookrightarrow}

\nc{\iso}{\overset{\sim}{\lra}}

\nc{\ssubset}{\underset{\not=}{\subset}}

\nc{\vac}{|0\rangle}

\nc{\Thm}[1]{Theorem~\ref{#1}}
\nc{\Prop}[1]{Proposition~\ref{#1}}
\nc{\Lem}[1]{Lemma~\ref{#1}}
\nc{\Cor}[1]{Corollary~\ref{#1}}
\nc{\Conj}[1]{Conjecture~\ref{#1}}
\nc{\Claim}[1]{Claim~\ref{#1}}
\nc{\Defn}[1]{Definition~\ref{#1}}
\nc{\Exa}[1]{Example~\ref{#1}}
\nc{\Rem}[1]{Remark~\ref{#1}}
\nc{\Note}[1]{Note~\ref{#1}}
\nc{\Quest}[1]{Question~\ref{#1}}
\nc{\Hyp}[1]{Hypoth\`ese~\ref{#1}}
\nen{thm}[1]{\label{#1}{\bf Theorem.\ } \em}{}
\nen{prop}[1]{\label{#1}{\bf Proposition.\ } \em}{}
\nen{lem}[1]{\label{#1}{\bf Lemma.\ } \em}{}
\nen{cor}[1]{\label{#1}{\bf Corollary.\ } \em}{}
\nen{conj}[1]{\label{#1}{\bf Conjecture.\ } \em}{}

\nen{claim}[1]{\label{#1}{\bf Claim.\ } \em}{}

\nen{defn}[1]{\label{#1}{\bf Definition.\ } }{}
\nen{exa}[1]{\label{#1}{\bf Example.\ } }{}
\nen{rem}[1]{\label{#1}{\em Remark.\ } }{}
\nen{note}[1]{\label{#1}{\em Note.\ } }{}
\nen{exer}[1]{\label{#1}{\em Exercise.\ } }{}
\nen{sket}[1]{\label{#1}{\em Sketch of proof.\ } }{}
\nen{quest}[1]{\label{#1}{\bf Question.\ } \em}{}

\nen{hyp}[1]{\label{#1}{\bf Hypoth\`ese.\ } \em}{}
\begin{document}


\title[]{On a generic Verma module at the critical level over
affine Lie superalgebras}

\author[Maria Gorelik]{{Maria Gorelik}
}
\address{
{\tt email: maria.gorelik@weizmann.ac.il} 
}

\begin{abstract}
We describe the structure of a Verma module with a generic highest weight 
at the critical level over a  symmetrizable affine Lie superalgebra
$\fhg\not=A(2k,2l)^{(4)}$.  We obtain the character formula for
a simple module with a generic highest weight 
at the critical level conjectured by V.~G.~Kac and D.~A.~Kazhdan.
\end{abstract}

\thanks{The author was partially supported by TMR Grant No. FMRX-CT97-0100}
\maketitle

\section{Introduction}
It is well-known that the representation theory of 
a complex  affine Lie algebra
changes drastically at the critical level.
In particular, Verma modules 
contain infinite number of singular vectors of imaginary degrees.
As it is shown in~\cite{Ku},
a Verma module with a ``generic'' highest weight 
at the critical level looks like a polynomial algebra
in a countable number of variables: submodules correspond
to the ideals in the polynomial algebra and the Jantzen filtration
corresponds to the adic filtration. Here ``genericity'' 
means that  all singular vectors
lie in the imaginary degrees.
As a consequence, J.-M.~Ku obtains the character formula conjectured
by ~V.~G.~Kac and D.~A.~Kazhdan in~\cite{kk}: 
$\ch L(\lambda)=e^{\lambda}\prod_{\alpha\in\hat{\Delta}^+_{re}} 
(1-e^{-\alpha})^{-1}$, where 
 $\lambda$ is a generic highest weight 
at the critical level and $L(\lambda)$ is the corresponding simple module.
In this paper we extend the results
of~\cite{Ku} to symmetrizable affine Lie superalgebras 
$\fhg\not=A(2k,2l)^{(4)}$ (see~\ref{in1}).

For non-twisted affine Lie algebras the
Kac-Kazhdan character formula was proven by different methods:
for $\hat\fsl(2)$ by M.~Wakimoto~\cite{Wk},
N.~Wallach~\cite{Wl}; for the affinizations
of classical algebras by T.~Hayashi~\cite{H} and 
R.~Goodman, N.~Wallach~\cite{GW}; for the affinization of
a general simple Lie algebra by 
B.~Feigin and E.~Frenkel~\cite{FF},~\cite{F} and recently
by T.~Arakawa~\cite{Ar}; in finite characteristic
by  O.~Mathieu~\cite{M}. For an arbitrary affine Lie algebras  
(including the twisted case) the formula was proven by
J.-M.~Ku~\cite{Ku} and recently reproven by M.~Szczesny~\cite{Scz}.

The approach of B.~Feigin, E.~Frenkel and M.~Szczesny
is based on the explicit realization of $L(\lambda)$:
they show that if $\lambda$ is a generic highest weight 
at the critical level then $L(\lambda)$ is isomorphic to a 
Wakimoto module, which is a representation
of $\fhg$ in a Fock module over some infinite-dimensional 
Heisenberg algebra; the construction of Wakimoto modules uses
a technique of vertex algebras. The Heisenberg algebra 
here corresponds to the set of {\em real} roots of $\fhg$.
The method of J.-M.~Ku is much more ``elementary'': it is based 
on a study of singular vectors in the Verma module $M(\lambda)$.
It can be interpreted (see~\ref{inH})
in terms of an infinite-dimensional 
Heisenberg algebra which corresponds to the set of {\em imaginary}
 roots of $\fhg$ 
(this Heisenberg algebra is a subalgebra of $\fhg$).
Our approach is close to one of J.-M.~Ku.

\subsection{Main result}\label{mainres}
Let $\fhg=\fhn_-\oplus\fhh\oplus\fhn$ be an affine Lie superalgebra
with a symmetrizable indecomposable
Cartan matrix (see~\ref{in1}).
Let $M(\lambda)$ be a Verma module of the highest weight $\lambda$
and $v_{\lambda}$ be its canonical generator, which we assume to
be even.

The Lie superalgebra $\fhn_-$ admits a triangular decomposition
$\fhn_-=\cN^-_-\oplus\cH_-\oplus\cN^+_-$, where $\cH_-$ consists
of the elements of imaginary weights (for a non-twisted case,
$\cH_-=\cL\fh\cap \fhn_-$ and $\cN^{\pm}_-=\cL\fn_{\pm}\cap \fhn_-$,
where $\cL$ stands for the loop space of a given subalgebra of
$\fg$). Set
$$\cS:=\cU(\cH_-).$$ 
Introduce the projections $\HC_{\pm}:\cU(\fhn_-)\to \cS$, where 
$$\Ker\HC_+=\cU(\fhn_-)\cN^-_-+\cN^+_-\cU(\fhn_-),\ \ 
\Ker\HC_-=\cU(\fhn_-)\cN^+_-+\cN^-_-\cU(\fhn_-).$$
For each $\lambda\in\fhh^*$  
define $\HC_{\pm}: M(\lambda)\to \cS$ via the natural
identification of $M(\lambda)$ with $\cU(\fhn_-)$;
these  projections
play a central role in our description of $M(\lambda)$.
The  projection $\HC_+$ appeared  in~\cite{Ku} and~\cite{Ch}. 

We call $v\in M(\lambda)$ {\em singular} 
if $v$ is a weight vector and
$v\in M(\lambda)^{\fhn}$ (these vectors are also called
primitive).
We say that {\em $\lambda$ is a critical weight} or 
{\em $\lambda$ has the critical level} if
$(\lambda+\hat{\rho},\delta)=0$ for an
imaginary root $\delta$. If $\lambda$ is a critical weight then
$M(\lambda)_{\lambda-\delta}$ contains a singular vector;
we say that $\lambda$ is a generic critical weight
($\lambda\in\Lambda_{crit}$) if
$M(\lambda)_{\lambda-\alpha}$ has no singular vectors
unless $\alpha$ is  proportional to $\delta$.

If $\lambda$ is a generic critical weight, 
the space $M(\lambda)^{\fhn}$ has a natural structure of an
associative (super)algebra. 
Indeed, the map $\psi\mapsto \psi(v_{\lambda})$
gives an embedding of $\End_{[\fhg,\fhg]}(M(\lambda))$
into the  space of singular vectors $M(\lambda)^{\fhn}$. 
The genericity condition on $\lambda$
means that this embedding is bijective and this provides
$M(\lambda)^{\fhn}$ with the algebra structure.

{\em In~\ref{thm01}--\ref{janf} we assume that $\fhg$
is not of type $A(2k,2l)^{(4)}$ and that
$\lambda$ is a generic critical weight.}

\subsubsection{}
\begin{thm}{thm01}
Let $\lambda$ be a generic critical weight. 
\begin{enumerate}
\item One has
$[M(\lambda):L(\lambda-s\delta)]=\dim M(\lambda)^{\fhn}_{\lambda-s\delta}$;
thus any submodule of $M(\lambda)$ is generated by singular vectors.

\item The restrictions of $\HC_+: M(\lambda)\to \cS$ 
and of $\HC_-: M(\lambda)\to \cS$ to the space
of singular vectors $M(\lambda)^{\fhn}$ give algebra isomorphisms
 $M(\lambda)^{\fhn}\iso\cS$, where the image of $v_{\lambda}$
is $1\in\cS$.

\item Any singular vector generates a submodule isomorphic to 
$M(\lambda-s\delta)$ for some $s\geq 0$.
\end{enumerate}
\end{thm}

Notice that $\cH_-$ is even and commutative so $\cS$ is the algebra of
polynomials in countably many variables.

\subsubsection{}\label{HN}
For a submodule $N$ of $M(\lambda)$ set 
$$H(N):=\HC_+(N^{\fhn})\subset \cS.$$
From~\Thm{thm01} we see that $H$ provides a  one-to-one correspondence
between the submodules of $M(\lambda)$ and $\ad\fhh$-invariant ideals
of $\cS$. 

If $\lambda$ is a generic critical weight then all simple
subquotients of $M(\lambda)$ are of the form $L(\lambda-s\delta)$;
note that 
$L(\lambda)\cong L(\lambda-s\delta)$ as $[\fhg,\fhg]$-modules
so $\ch L(\lambda-s\delta)=e^{-s\delta}\ch L(\lambda)$.
As a result,
the characters of $N$ and of $H(N)$ are connected by
the following formula:
\begin{equation}\label{cor01}
\ch N=\ch L(\lambda)\cdot\ch H(N).
\end{equation}

Applying this formula to $N=M(\lambda)$ we get the  
Kac-Kazhdan  character formula:
$$\ch M(\lambda)=\ch L(\lambda)\cdot\ch \cS, $$ 
that is
$$\ch L(\lambda)=e^{\lambda}\prod_{\alpha\in\hat{\Delta}^+_{re;\ol{0}}} 
(1-e^{-\alpha})^{-1}
\prod_{\alpha\in\hat{\Delta}^+_{re;\ol{1}}} (1+e^{-\alpha}).$$

Note that $\ch L(\lambda)=e^{\lambda}\ch \cU(\cN^+_-)\ch\cU(\cN^-_-)$.

\subsubsection{Jantzen filtration}\label{janf}
Recall that $\cH_-$ is commutative so $\cS=\cU(\cH_-)$ 
is the symmetric algebra: $\cS=\sum_{j=0}^{\infty}\cS^j$.
The spaces $\cS^{\geq k}:=\sum_{j=k}^{\infty}\cS^j$
form the adic  filtration on $\cS$.

\begin{thm}{thm02} If $\lambda$ is a generic critical weight then
$H$ maps the Jantzen filtration
$\{M(\lambda)^k\}$ to the adic filtration on $\cS$,
i.e. $H(M(\lambda)^k)=\cS^{\geq k}$.
\end{thm}

As a result, the Jantzen filtration
coincides with the cosocle filtration.

Combining Theorems~\ref{thm01}, \ref{thm02} we obtain the

\begin{cor}{}
$M(\lambda)^k$ is generated by the singular vectors
$\HC_+^{-1}(\cS^k)$.
\end{cor}

\subsection{Connection with representations of Heisenberg algebra}
\label{inH}
Let $\cH$ (resp., $\cH_-$) be the sum of  positive (resp., negative)
imaginary root spaces of
$\fhg$. Set $\fl=\cH_-\oplus\mathbb{C}K\oplus\mathbb{C}D
\oplus\cH$ and let $V^k$ be the induced module
$V^k=\Ind_{\mathbb{C}K\oplus\mathbb{C}D
\oplus\cH}^{\fl}\mathbb{C}_k$ from the one-dimensional module
$\mathbb{C}_k$ with trivial action of
$\cH+\mathbb{C}D$ and $K=k\in\mathbb{C}$. Identify $V^k$
with $\cS=\cU(\cH_-)$ and say that $v\in V^k$ is singular
if $\cH v=0$ and $v$ is an $D$-eigenvector.

Let $\fhg\not=A(2k,2l)^{(4)}$. Then $\cH_-\oplus\mathbb{C}K\oplus\cH$
is a countably dimensional Heisenberg algebra and
$V^k$ is irreducible if and only if $k\not=0$. In $V^0$
any $D$-eigenvector is singular;
the Jantzen filtration of $V^0$ identifies
with the adic filtration of $\cS$.
Theorems~\ref{thm01},\ref{thm02} can be reformulated in the following 
way.

\subsubsection{}\label{reH}
{\em $M(\lambda)$ with a generic highest weight 
looks like $V^k$ for $k:=(\lambda+\hat{\rho},\delta)$. More precisely, 
$\HC_+: M(\lambda)\to \cS$ maps singular vectors to the singular ones
and induces a bijection  between the submodules
of $M(\lambda)$ and of $V^k$.
This bijection is compatible with the Jantzen filtrations.}

\subsubsection{}
In~\cite{GS} the result similar to~\ref{reH} 
is proven for a non-symmetrizable
affine Lie superalgebra $q(n)^{(2)}$.
In this case $V^k$ is reducible for all $k\in\mathbb{C}$
and $M(\lambda)$ is reducible
for all $\lambda\in\fhh^*$.

\subsubsection{}
We believe that a result similar to~\ref{reH} holds for 
$\fhg=A(2k,2l)^{(4)}$. A difficulty in this case is that 
$\cH_-$ is not commutative and it is not true that any submodule of
$V^k$ (and of $M(\lambda)$) is generated by
singular vectors.

\subsection{Affine Lie superalgebras}\label{in1}
Recall that affine Lie algebras are  
finite growth Kac-Moody  algebras. 
The Cartan matrix of an
affine Lie algebra is symmetrizable~\cite{Kacbook}.
An affine Lie algebra can be described 
in terms of a finite-dimensional semisimple Lie algebra and its 
finite order automorphism (see, for example,~\cite{Kacbook},
Ch. VI-VIII).  

The superalgebra generalization of Kac-Moody algebras
was introduced in~\cite{Kadv}; a detailed treatment of this notion
can be found in~\cite{Wa}. 
Call a Kac-Moody superalgebra
{\em affine} if it has a finite growth and {\em 
symmetrizable} if it has a symmetrizable Cartan matrix.
In~\cite{vdL},\cite{S}
the affine symmetrizable  Lie superalgebras 
are described
 in terms of  finite-dimensional Kac-Moody  superalgebras 
and  their finite order automorphisms; we recall this construction
in~\ref{affsym}.

A symmetrizable affine Lie superalgebra has 
a Casimir element. As a consequence, a Verma module 
$M(\lambda)$ is irreducible unless $\lambda$ belongs to 
the union of countably many hyperplanes. Among these hyperplanes
one is rather special: $M(\lambda)$ has an infinite length
if $\lambda$ lies on this hyperplane. This hyperplane is the set of
critical weights.

Affine  Lie superalgebras were classified in a recent paper~\cite{HS}. 
It turns out that non-symmetrizable affine  Lie superalgebras 
consist of 4 series. One of this series is $q(n)^{(2)}$;
these algebras are twisted affinizations of ``strange'' Lie superalgebras
 $q(n)$. As it is mentioned above, a Verma module over $q(n)^{(2)}$
is always reducible. It is an interesting question to study
Verma modules of generic highest weight for other series.

\subsection{Outline of the proof}
In~\ref{out1},~\ref{out2} below we outline the proof of
theorems~\ref{thm01},~\ref{thm02}. Retain notation of~\ref{mainres}.

\subsubsection{}\label{out1}
In~\ref{reform}
we show that theorems~\ref{thm01},~\ref{thm02} can be easily deduced
from the following assertions for the Verma module
of a generic critical weight $\lambda$: 

(A) the restrictions of
$\HC_{\pm}$ to the set of singular vectors of $M(\lambda)^{k}$ 
contain $\cS^{\geq k}$,

(B) $[M(\lambda)^k:L(\lambda-s\delta)]
=\dim\cS^{\geq k}_{-s\delta}$ for all $k,s$.

We prove (A), (B) as follows. 
Let $A$ be the local algebra
$A=\mathbb{C}[x]_{(x)}$ and $M(\lambda+x\xi)$ be a Verma module 
over $\fhg_A=\fhg\otimes A$.
 Recall that the Jantzen filtration 
$M(\lambda)^{k}$ is obtained by the evaluation of the Jantzen filtration 
$\cF^k(M(\lambda+x\xi))$ at $x=0$.

Let $\sigma$ be the natural anti-involution of $\fhg$ 
(it interchanges the canonical generators $e_i, f_i$
and acts by $\id$ to $\fhh$). Set $\cN^{\pm}=\sigma(\cN^{\mp}_-)$.

For $\lambda$  generic 
(this means that either $\lambda$ is a generic critical weight 
or $M(\lambda)$ is simple) we construct linear
maps $\HC_{\pm,A}^{-1}:\cS_A\to M(\lambda+x\xi)^{\cN^{\pm}}$.
Denote by $\HC_{\pm}^{-1}:\cS\to M(\lambda)$ the reductions
of $\HC_{\pm,A}^{-1}$ modulo $(x)$.
We show that the images of $\HC_{\pm}^{-1}$
lie in $M(\lambda)^{\fhn}$ if $\lambda$ is a generic critical weight.

Let $(-,-)$ be the Shapovalov form on $M(\lambda+x\xi)$.
One readily sees that $M(\lambda+x\xi)^{\cN^{+}}$
is orthogonal to $\Ker\HC_-$ with respect to
this form. As a result, $\HC_{+,A}^{-1}(z)\in\cF^k(M(\lambda+x\xi))$
if $(\HC_{+,A}^{-1}(z),\HC_{-,A}^{-1}(\cS))\subset (x)^k$.
We prove the last inclusion for
$k$ equal to the degree of $z$ by induction on this degree; 
the proof is based on
the fact that $\HC_-^{-1}(z')$ is singular for any $z'$. 
This is the point where
we use simultaneously both projections $\HC_{\pm}$
in order to prove a certain assertion for $\HC_+$.
This gives (A) for $\HC_+$; the proof for $\HC_-$
is similar.

The proof of (B) is based on (A) and the formula:
$$\sum_{k\geq 1}\ch M(\lambda)^k =\frac{\ch M(\lambda)}{\ch \cS}
\sum_{k\geq 1}\ch \cS^{\geq k}$$
which easily follows from the Jantzen sum formula.
We deduce (B) from this formula by induction on weight.

\subsubsection{}\label{out2}
One of the key points of the proof is the construction of
$\HC_{+,A}^{-1}(z)$ for $z\in\cH_-$ and 
the proof that $\HC_+(z)$ is singular if $\lambda$
is a generic critical weight.

Let $\lambda$ be generic.
The construction of $\HC_{+,A}^{-1}(z)$ for $z=a(-m)\in\cH_-$
($m>0$) follows the construction of singular vectors
in $M(\lambda)$ given in~\cite{Ku}; we outline this construction
in~\ref{outv+}.

The proof that $\HC_+^{-1}(z)$ is singular
goes as follows. It is easy to see that for generic $\lambda$ a vector
$v\in M(\lambda)_{\lambda-m\delta}$
is singular if $(\cN^++\cH)v=0$. By the construction
$v:=\HC_+^{-1}(z)$
is annihilated by $\cN^+$ and by the elements of $\cH_{s\delta}$
for $s\not=m$. From~\ref{symN} it is easy to deduce
that for $u\in\cH_m$ one has
$$uv=(u|a(-m))(mK+c)v_{\lambda},$$
where $c$ is a scalar which does not depend on $\lambda$.
Therefore $v$ is singular if $\lambda$ has level 
$-c/m$. Recall that $M(\lambda)$ is simple if
the level is not critical so $v$ is not singular 
in this case. Hence $v$ is singular if the level of $\lambda$
is critical.

\subsection{Acknowledgment}
I am very grateful to M.~Duflo for fruitful discussions. 
A part of this work was done during  my stay at IHES and at MPIM.
I am grateful to these institutions for stimulating atmosphere and excellent
working conditions.

\section{Preliminaries and notation}
\label{prelim}
Our base field is $\mathbb{C}$; we denote by $\mathbb{Z}_{\geq 0}$
the set of non-negative integers. If $V$ is a superspace, 
we denote by $p(v)$ the parity of a vector $v\in V$.
For a Lie superalgebra $\fg=\fg_{\ol{0}}\oplus\fg_{\ol{1}}$ 
we denote by $\Ug$ its universal enveloping
algebra and  by $\Sg$ its symmetric algebra.
We view $\Ug$ as a $\fg$-module with respect to the adjoint action.

If $\fp$ is a Lie algebra, $N$ is a 
$\fp$-module and $N'$ is a subspace of
$N$, denote by $\fp N'$ the vector space spanned by
$xv$ where $x\in \fp,\ v\in N'$.

\subsection{Index of notations}
\label{index}
Symbols used frequently are given below under the section number
where they are first defined. 
$$\begin{array}{llccll}
\ref{trfhg}&\tilde{\cH},\ \tilde{\cN^{\pm}}, & & & 
\ref{cHN} & \cH,\ \cN^{\pm},\ \cH_-,\ \cN^{\pm}_-,\\
\ref{cS}& \cS,\ \cS^k,\ \cS^{\geq k}
& & &  \ref{HC--}& \HC_+,\ \HC_-,\ \breve{M}
  \\
\ref{Omega} & \Omega(M),\ \wt u & & & 
\ref{defgn} & \Lambda,\ \Lambda_{crit}.
\end{array}$$

\subsection{Triangular decompositions of superalgebras}
\label{triang}
A triangular decomposition of a Lie superalgebra $\fg$
can be constructed as follows (see~\cite{PS}).
A Cartan algebra is a nilpotent subalgebra which
coincides with its normalizer. 
Fix a Cartan algebra $\fh$. Then $\fg$ has a generalized root
decomposition
$\fg:=\fh + \oplus_{\alpha\in\Delta}{\fg}_{\alpha}$,
where $\Delta$ is a subset of ${\fh}^*$ and
$$\fg_{\alpha}=\{x \in \fg|(ad (h)-\alpha (h))^{dim \fg}(x)=0\}.$$

In this paper all Cartan algebras are pure even
and any root space ${\fg}_{\alpha}$ is
either odd or even. This allows one to define the parity on the set of
roots $\Delta$: we
denote by $\Delta_0$ (resp. $\Delta_1$) the 
set of non-zero weights of $\fg_{\ol{0}}$ (resp., $\fg_{\ol{1}}$) with
respect to $\fh$; one has $\Delta=\Delta_0\coprod\Delta_1$.

Now fix $h\in\fh_0$ satisfying $\Ree \alpha(h)\not=0$ for $\alpha(h)\not=0$
for any $\alpha\in\Delta$ (i.e. $\alpha(h)$ 
is not a non-zero imaginary number). Set 
$$\begin{array}{ll}
\Delta^+:=\{\alpha\in\Delta|\ \Ree \alpha(h)>0\},\\
\fn:=\sum_{\alpha\in\Delta^+}\fg_{\alpha}.
\end{array}$$
Define $\Delta^-$ and $\fn_-$ similarly. Then 
$\fg=\fn_-\oplus\fh\oplus\fn$
is a triangular decomposition.

\subsection{Affine Lie superalgebras with symmetrizable Cartan matrices}
\label{affsym}
According to~\cite{vdL}, any Kac-Moody  superalgebra,
which has finite Gelfand-Kirillov
dimension and a symmetrizable Cartan matrix, can be described in terms
of the loop algebra of a   finite dimensional Kac-Moody
Lie superalgebra and its automorphism of a finite order.
The non-twisted affine Lie superalgebras correspond
 to the trivial automorphism. 
The twisted affine  Lie superalgebras correspond to the Dynkin  diagrams
 of types $X_N=A(2k,2l-1)$, $A(2k-1,2l-1)\ \bigl((k,l)\not=(1,1)\bigr)$, 
$A(2k,2l)$, $C(l+1)$, $D(k+1,l)$, $G_3$ 
and automorphisms of order $2,2,4,2,2,2$
respectively. We briefly recall the construction below.

\subsubsection{}\label{XN}
Let $X_N$ be a  Dynkin  diagram of a Kac-Moody
finite dimensional Lie superalgebra $\fg=\fn_-\oplus\fh\oplus\fn$
 with the triangular decomposition determined by this Dynkin diagram.
We may (and will) assume that $X_N$ is connected.
Then $\fg$ is either simple or of the type $\fgl(n|n)=A(n-1,n-1)$.
The even part $\fg_{\ol 0}$ is reductive. The algebra $\fg$
admits a non-degenerate even invariant bilinear form $(-|-)'$;
if $\fg$ is of the type $A(m,n)$ 
we take the form $(a|b)'=\str ab$.

\subsubsection{}\label{XN1}
The affine Lie superalgebra $\fhg'$ 
corresponding to  the Dynkin diagram $X_N^{(1)}$  can be described as
follows:
$\fhg'=\displaystyle\oplus_{m\in\mathbb{Z}}\fhg'_m$,
where $\fhg'_0=\fg\oplus \mathbb{C} K'
\oplus \mathbb{C} D'$ and, for $m\not=0$, $\fhg'_m:=\fg t^m$ if
$\fg\not=\fgl(n|n)$, $\fhg'_m:=\fpgl(n|n)t^m$ for $\fg=\fgl(n|n)$.
The structure of superspace on $\fhg'$
is given by $p(t)=p(K')=p(D')=0$ and the commutational relations are
$$\begin{array}{l}
\ [a(m),b(k)]=[a,b](m+k)+m\delta_{m,-k}(a|b)'K',\ \ [\fhg,K']=0,\ \ 
[D',a(m)]=ma(m),
\end{array}$$
where $a(m):=at^m$ for $a\in\fg, m\in\mathbb{Z}$.
If $\fg=\fgl(m,m)$ then the term  $[a,b]$
should be substituted by its image in $\fpgl(n|n)$
if $m+k\not=0$ and $mk=0$.

The form $(-|-)'$ 
can be extended from $\fg\subset \fhg'_0$ to $\fhg'$ by setting
$(K'|D')'=0,\ (K'|D')'=1$ and
$(xt^m|yt^n)'=\delta_{m,n}(x|y)',\ (\mathbb{C}K'+\mathbb{C}D'| xt^m)'=0$
for any $x,y\in\fg$.

\subsubsection{}\label{XNr}
Let $\epsilon$ be an automorphism of $\fg$ of a finite order
$r$ ($r=1,2,3,4$) such that $\epsilon(\fh)=\fh$.
Extend $\epsilon$ to $\fhg'$ 
by putting $\epsilon(t)=\exp(-2\pi i/r)t$, that is 
$$\epsilon(at^m )=\exp(-2\pi mi/r)\epsilon(a)t^m,\ \ 
\epsilon(K')=K',\ \ \epsilon(D')=D'.$$
By~\cite{Kacbook},\cite{vdL} the affine Lie superalgebra $\fhg$ 
 corresponding to the Dynkin diagram $X_N^{(r)}$ 
is isomorphic to the subalgebra of invariants
$(\fhg')^{\epsilon}$. We identify this subalgebras:
$$\fhg:=(\fhg')^{\epsilon}.$$

Clearly, $\fhg$ inherits the $\mathbb{Z}$-grading:
$\fhg=\displaystyle\oplus_{m\in\mathbb{Z}}\fhg_m$, where
$\fhg_m=(\fhg'_m)^{\epsilon}$.

Introduce $K:=rK', D:=D'$ and the form $(-|-):=(-|-)'/r$.
The restriction of $(-|-)$ gives an invariant form on $\fhg$.
We identify $\fg^{\epsilon}\oplus\mathbb{C}K\oplus\mathbb{C}D$
with the zero homogeneous component $\fhg_0$.

\subsubsection{}\label{trfg}
Recall (\cite{Kacbook}, 8,1) that any automorphism of 
a finite order of 
the reductive algebra $\fg_{\ol{0}}$ has an invariant regular element.
Without loss of generality we choose the triangular decomposition
$\fg_{\ol{0}}=\fn_{-,\ol{0}}\oplus\fh\oplus\fn_{\ol{0}}$
determined by this regular element.

Then $\fh^{\epsilon}$
is a Cartan subalgebra of $\fg^{\epsilon}$ and 
the centralizer of $\fh^{\epsilon}$ in $\fg_{\ol{0}}$ coincides with $\fh$. 
We have a triangular decomposition
$\fg_{\ol{0}}^{\epsilon}=\fn_{-,\ol{0}}^{\epsilon}\oplus
\fh^{\epsilon}\oplus\fn_{\ol{0}}^{\epsilon}$.
Recall that all triangular decompositions of $\fg_{\ol{0}}^{\epsilon}$
are conjugate by inner automorphisms of $\fg_{\ol{0}}^{\epsilon}$.
Since any inner automorphism of $\fg_{\ol{0}}^{\epsilon}$
can be extended to $\fhg$, any triangular 
decomposition
$$\fhg=\fhn_-\oplus\fhh\oplus\fhn$$
is conjugate by $\Aut(\fhg)$ to a triangular decomposition 
compatible with $\fg_{\ol{0}}^{\epsilon}=\fn_{-,\ol{0}}^{\epsilon}\oplus
\fh^{\epsilon}\oplus\fn_{\ol{0}}^{\epsilon}$.
Hence without loss of generality we may (and will) 
assume that
$$\fhh=\fh^{\epsilon}\oplus\mathbb{C}K\oplus\mathbb{C}D,\ \ 
\fn_{\ol{0}}^{\epsilon}\subset \fhn.$$

By above $\fh$ is
the centralizer of $\fh^{\epsilon}$ in $\fg_{\ol{0}}$.
If $\fhg$ is not of the type $A(2k,2l)^{(4)}=\fsl(2k+1,2l+1)^{(4)}$  
the centralizer of $\fh^{\epsilon}$ in $\fg_{\ol{1}}$ is zero.
For $A(2k,2l)^{(4)}$ the centralizer of $\fh^{\epsilon}$ in
$\fg_{\ol{1}}$ has dimension two. Hence

\begin{equation}\label{strcH}\begin{array}{l}
\fhg^{\fh^{\epsilon}}=\fhh+\sum_{k\not=0} (\fh t^k)^{\epsilon}, \ 
\text{ for } \fhg\not=A(2k,2l)^{(4)},\\
\fhg_{\ol{0}}^{\fh^{\epsilon}}=\fhh+\sum_{k\not=0} (\fh t^k)^{\epsilon}, \ 
\text{ for } \fhg=A(2k,2l)^{(4)}.
\end{array}
\end{equation}

\subsection{Roots}
Denote by $\Delta^+$ (resp.,  $\hat{\Delta}^+$)
the multiset of positive roots of $\fg$ (resp., of $\fhg$), and
by $\hat{\Delta}^+_0$ (resp., $\hat{\Delta}^+_1$) the multiset of
even (resp., odd) positive roots. Remark that $\Delta^+$ is a usual set
(all root spaces in $\fg$ are one-dimensional).

The bilinear form $(-|-)$ induces 
a non-degenerate bilinear form $(-,-)$ on $\fhh^*$.
Let $h_{\alpha}$ be the image of $\alpha$ under
the isomorphism $\fhh^*\to\fhh$ induced by the bilinear form
(that is $\mu(h_{\alpha})=(\alpha,\mu)$).

\subsubsection{}
Real roots are defined via the action of Weyl group 
(see~\cite{Wa} for the definition)
on $\hat{\Delta}^+$; the roots is imaginary if it is not real.
For a symmetrizable affine Lie superalgebras 
a root $\alpha\in \hat{\Delta}^+$ is imaginary if $(\alpha,\beta)=0$
for all $\beta\in\hat{\Delta}$.
Denote  by $\hat{\Delta}^+_{re}$ (resp., $\hat{\Delta}^+_{im}$)
the multiset of real (resp., imaginary) positive roots;
one has $\hat{\Delta}^+=\hat{\Delta}^+_{re}\coprod\hat{\Delta}^+_{im}$.
All imaginary roots are proportional to a certain element  
$\delta\in\fhh^*$:
if $\fhg_{\alpha}\subset\fhg_s$ and $\alpha$
is imaginary then $\alpha=s\delta$. Notice that $h_{\delta}=K$,
since $h_{\delta}$ is central and that $(D|h_{\delta})=D(\delta)=1$.

Apart of the case
$G(3)^{(2)}$ one has $\hat{\Delta}^+_{im}=\mathbb{Z}\delta$;
for $G(3)^{(2)}$ one has $\hat{\Delta}^+_{im}=\mathbb{Z}(2\delta)$,
i.e. $\delta$ is not a root.  All imaginary roots are even
if $\fhg$ is not of the type $A(2k,2l)^{(4)}$.

For a non-twisted case 
$\hat{\Delta}^+_{re}=\mathbb{Z}\delta+\Delta$, 
$\fg_{\alpha+m\delta}=\fg_{\alpha}t^m$ and
$\fhg_{m\delta}=\fh t^m$ for $m\not=0$.

\subsubsection{}
Let $\hat{\pi}$ be the set of simple roots for $\hat{\Delta}^+$.
Let $\hrho\in\fhh$ be such that $(\hrho,\alpha)=\frac{1}{2}(\alpha,\alpha)$
for any $\alpha\in\hat{\pi}$.
For a weight $\lambda\in\fhh^*$ the value $\lambda(K)$
is called a {\em level}. The weights of a given level form a hyperplane
in $\fhh^*$. Put $h^{\vee}:=\hrho(K)$.
The level $-h^{\vee}$ is called {\em critical}:
$\lambda$ has the critical level iff $(\lambda+\hrho,\delta)=0$.

\subsubsection{}
Set $\hat{Q}^+:=\sum_{\alpha\in\hat{\Delta}^+}\mathbb{Z}_{\geq 0}\alpha$.
Define a partial ordering on $\fhh^*$ by setting
$\mu\geq\mu'$ if $\mu-\mu'\in\hat{Q}^+$.

\subsection{Algebras $\cH,\cN^{\pm}$}
\label{trfhg}
One has $\fhg=[\fhg,\fhg]\oplus\mathbb{C}D$.
Retain notation of~\ref{XNr} and
set
$$\fh':=\fhh\cap [\fhg,\fhg].$$
One has $\fh'=\fh^{\epsilon}+\mathbb{C}K$.
Observe that the decomposition $[\fhg,\fhg]=\fhn_-\oplus \fh'
\oplus\fhn$ is not triangular in a sense of~\ref{triang},
 since the centralizer of $\fh'$ is
$$\tilde{\cH}:=\fh'+\sum_{j\not=0} \fhg_{j\delta}.$$
Fix $h'\in\fh'$ in such a way that 
$\Ree\alpha(h')>0$ if $\alpha$ is a weight of $\fhn\cap\fhg_0$
and that $\Ree\alpha(h')\not=0$ for any $\alpha\in\hat{\Delta}_{re}$.
Fix a triangular decomposition 
$$[\fhg,\fhg]=\tilde{\cN}^+\oplus\tilde{\cH}\oplus\tilde{\cN}^-$$
defined via $h'$ as in~\ref{triang}. Then
$\fhn\cap\fhg_0\subset\tilde{\cN}^+,\ \fhn_-\cap\fhg_0\subset\tilde{\cN}^-$.
If $\fhg$ is a non-twisted affinization one has
$\tilde{\cN}^{\pm}=\sum_{m\in\mathbb{Z}}\fn_{\pm}t^m$.

Let $\fhg$ is not of the type $A(2k,2l)^{(4)}$. 
By~(\ref{strcH}), $\tilde{\cH}$ is spanned  elements of the form $ht^k,
h\in\fh$ and by $K$.
In particular, $\tilde{\cH}$ is isomorphic to the direct product of
countably dimensional Heisenberg algebra 
$\mathbb{C}K+\sum_{j\not=0} \fhg_{j\delta}$
and the abelian Lie algebra $\fh^{\epsilon}$.

\subsection{Automorphism of $\fhg$}\label{aut}
In~\Lem{lemsym} we will show that for each $\fhg$  
there exist a commuting pair $\epsilon,
\phi\in\Aut(\fg)$ such that $\fhg=(\fhg')^{\epsilon}$ 
in a sense of~\ref{XNr} and $\phi|_{\fh}=-\id$.

Since $\phi|_{\fh}=-\id$, the  bilinear form $(-,-)$
is $\phi$-invariant and this allows us to define an automorphism 
$\hat{\phi}$ of $\fhg'$ by the formulas: 
$$\hat{\phi}(xt^k)=\phi(x)t^{k},\ \ \hat{\phi}(K)=K,\ \ \hat{\phi}(D)=D.$$
Since $\phi$ commutes with $\epsilon$,
the automorphism $\hat{\phi}$ stabilizes $\fhg$.
View the restriction of $\hat{\phi}$ to $\fhg$
as an automorphism of $\fhg$. Set
$$(\mathbb{C}K+\mathbb{C}D)^{\perp}:=\{h\in\fhh|\ 
(h,K)=(h,D)=0\}.$$

\subsubsection{}
\begin{cor}{symN}
$\hat{\phi}\in\Aut(\fhg)$ has the following properties
$$\begin{array}{ll}
(i) & \hat{\phi}(K)=K,\ \  \hat{\phi}(D)=D,\ \
\hat{\phi}|_{(\mathbb{C}K+\mathbb{C}D)^{\perp}}=-\id;\\
(ii) & \text{ for } \fhg\not=A(2k,2l)^{(4)}\ \ 
\hat{\phi}|_{\tilde{\cH}_j}=-\id,\ \text{ for }j\not=0;\\
(ii') & \text{ for } \fhg=A(2k,2l)^{(4)}\ \ 
\hat{\phi}|_{\tilde{\cH}_{2j}}=-\id,\ \text{ for }j\not=0,\ \
\hat{\phi}|_{\tilde{\cH}_{4k\pm 1}}=\pm\sqrt{-1}\id;\\
(iii) & \hat{\phi}(\tilde{\cN}^+_k)=\tilde{\cN}^-_k,\ \ 
\hat{\phi}(\tilde{\cN}^-_k)=\tilde{\cN}^+_k.
\end{array}$$
\end{cor}
\begin{proof}
Recall that $(\mathbb{C}K+\mathbb{C}D)^{\perp}=\fh^{\epsilon}$
and that for $j\not=0$ an element of $\tilde{\cH}_j$ 
(resp., $\tilde{\cH}_{2j}$)
takes  form $ht^j,\ h\in\fh$ for $\fhg\not=A(2k,2l)^{(4)}$
(resp., for $\fhg=A(2k,2l)^{(4)}$).
Hence $\hat{\phi}$ satisfies (i), (ii)
and the condition $\hat{\phi}|_{\tilde{\cH}_{2j}}=-\id$
for $\fhg=A(2k,2l)^{(4)}$. The property (iii) follows from (i).
For $\fg=A(2k,2l)^{(4)}$ the
 remaining condition
$\hat{\phi}|_{\tilde{\cH}_{4k\pm 1}}=\pm\sqrt{-1}\id$
easily follows from the fact that the spaces $\tilde{\cH}_{4k\pm 1}$
are one-dimensional.
\end{proof}

\subsubsection{}
\begin{lem}{lemsym}
There exists a commuting pair $\epsilon,
\phi\in\Aut(\fg)$ such that $\fhg=(\fhg')^{\epsilon}$ in a sense 
of~\ref{XNr} and $\phi|_{\fh}=-\id$.
\end{lem}
\begin{proof}
Let $\alpha_1,\ldots,\alpha_s$ be the set of simple roots of
$\fg$; for each $i$ choose $e_i\in\fg_{\alpha_i},\ f_i\in\fg_{-\alpha_i}$ 
such that $(e_i|f_i)=1$. For any collection of non-zero scalars
$\{c_i\}_{i\in I}$ define a map $\phi(C): \fg\to\fg$
by $\phi(e_i)=c_if_i,\ \phi(f_i)=(-1)^{p(e_i)}c_i^{-1}e_i$
and $\phi(h)=-h$ for $h\in\fh$.
Then $\phi(C)\in \Aut(\fg)$.

Let $\fg$ be a simple finite-dimensional Lie algebra.
Then we can choose $\epsilon$ induced by an automorphism
of the Dynkin diagram, that is $\epsilon(e_i)=e_{\sigma(i)}; 
\epsilon(f_i)=f_{\sigma(i)}$, where $\sigma$ is an automorphism
of the Dynkin diagram of $\fg$. Take $\phi$  
which acts by $-\id$ on $\fh$ and interchanges $e_i$ with $f_i$.
Clearly, $\phi$ commutes with $\epsilon$.

Let $\sigma$ be an automorphism of the Dynkin diagram of
a finite-dimensional Kac-Moody superalgebra. 
Let $D=\{d_i\}$ be the collection of elements
in $\mathbb{C}^*$. Denote by $\sigma(D)$ the automorphism
of $\fg$ given by the formulas $\sigma(D)(e_i)=d_ie_{\sigma(i)},\ 
\sigma(D)(f_i)=d_i^{-1}f_{\sigma(i)}$.
Suppose that $\sigma(D)^2=\id$. It is easy to check that $\sigma(D)$
commutes with $\phi(C)$ for the following collection
$C$: $c_i=d_i$ if $i\leq \sigma(i)$, $c_i=d_{\sigma(i)}^{-1}$
if $i>\sigma(i)$.

Hence any involution of the form $\sigma(D)$ commutes 
with  $\phi(C)$ for a suitable $C$.

Consider the case when
 $\fhg$ is a Lie superalgebra and $\fhg\not=A(2k,2l)^{(4)}$.
By~\cite{vdL}, Table 4, in this case $\epsilon^2=\id$ and
$\epsilon$ is of the form $\sigma(D)$ ($\sigma$ is identity for $G(3)$
and is an involution in other cases). By above,
$\phi(C)$ commutes with $\epsilon$ for a suitable collection $C$.

Consider the remaining case  $\fhg=A(2k,2l)^{(4)}$.
By~\cite{S} we can choose  $\epsilon$ to be $x\mapsto -x^{st}$, where
$st$ stands for
 the super-transposition in the matrix Lie superalgebra $\fg=A(2k,2l)$.
Then the restriction of $\epsilon$ to $\fg_{\ol{0}}$ is
$x\mapsto -x^t$, where $t$ stands for the transposition.
Hence  the restriction of $\epsilon$
to $\fg_{\ol{0}}$ is an involution.

Recall that $\epsilon$ leaves invariant a regular element of $\fh$
and fix the triangular decomposition
$\fg_{\ol{0}}=\fn_{\ol{0},-}\oplus \fh\oplus\fn_{\ol{0}}$
defined by this regular element. Then $\epsilon$ stabilizes $\fh$
and $\Delta^+_{\ol{0}}$.
Fix a Dynkin diagram
for $\fg$ containing a unique odd node (for $i=2k+1$)
and a triangular decomposition
which corresponds to this Dynkin diagram and is compatible
with the above triangular decomposition of $\fg_{\ol{0}}$.
Choose $\{e_i,f_i\}_{i\in I}$ as above (for this triangular decomposition).

Recall that $\fg_{\ol{0}}=\fsl(2k+1)\times\fsl(2l+1)\times\mathbb{C}$.
Note that $\{e_i,f_i\}_{i\in I, i\not=2k+1}$ generate 
$\fsl(2k+1)\times\fsl(2l+1)$. By above,
the restriction of $\epsilon$ to $\fsl(2k+1)\times\fsl(2l+1)$
is an involution which stabilizes $\fh$
and $\Delta^+_{\ol{0}}$. Therefore this
restriction takes form $\sigma(D)$, where $\sigma$ is an involution
of the Dynkin diagram $A(2k)\cup A(2l)$. By above,
for a suitable collection $C'=\{c_i\}_{i\not=2k+1}$ the automorphism
$\phi(C')\in\Aut(\fsl(2k+1)\times\fsl(2l+1))$ 
commutes with the restriction $\epsilon$  to $\fsl(2k+1)\times\fsl(2l+1)$.

For $\psi,\varphi\in\Aut(\fg)$ set 
$$[\psi,\varphi]:=\psi\circ\varphi\circ\psi^{-1}\circ\varphi^{-1}.$$

Take any $c_{2k+1}\in\mathbb{C}^*$ and let $\phi(C)\in\Aut(\fg)$
be the corresponding extension of $\phi(C')$.
Then  the restrictions of $\phi(C)$ and of $\epsilon$
to $\fsl(2k+1)\times\fsl(2l+1)$ commute. Moreover,
since $\epsilon$ stabilizes $\fh$ and $\phi(C)|_{\fh}=-\id$,
the restrictions of $\phi(C)$ and of $\epsilon$ to $\fh$
also commute. One has $\fg_{\ol{0}}=\fsl(2k+1)\times\fsl(2l+1)+\fh$
so $[\phi(C),\epsilon]|_{\fg_{\ol{0}}}=\id$.

Recall that 
${\fg}_{\ol{1}}=W\oplus W'$, where $W,W'$ are irreducible non-isomorphic
$\fg_{\ol{0}}$-modules. Any automorphism of $\fg$ which acts by $\id$
on $\fg_{\ol{0}}$ is equal to 
 $\psi(b)$  for $b\in\mathbb{C}^*$, where
$$\psi(b)|_{\fg_{\ol{0}}}=\id, \ \psi(b)|_{W}=b\id,\ 
 \psi(b)|_{W'}=b^{-1}\id.$$
By above, $[\phi(C),\epsilon]=\psi(a)$
for some $a$. Then for any $b\in\mathbb{C}^*$ one has
$$\begin{array}{rl}
[\psi(b)\circ\phi(C),\epsilon] &=\psi(b)\circ\psi(a)\circ
\psi(b)^{-1}\circ [\psi(b),\epsilon]=
\psi(a)\circ [\psi(b),\epsilon]\end{array}$$
Clearly, $[\psi(b),\epsilon]$ acts by $\id$ on $\fg_{\ol{0}}$. Moreover, 
$\epsilon$ interchanges $W$ with $W'$ so 
$[\psi(b),\epsilon]=\psi(b^2)$.
Now taking $b$ such that $ab^2=1$ we get
$[\psi(b)\circ\phi(C),\epsilon]=\psi(1)=\id$ that is
$\phi:=\psi(b)\circ\phi(C)$ commutes with $\epsilon$.
Since $\phi$ acts by $-\id$ on $\fh$ it satisfies
the required conditions.
\end{proof}

\subsection{}
\begin{lem}{pairing}
For any $x\in\tilde{\cH}_m, y\in \tilde{\cH}_{-m}$ one has
$$[x,y]=m(x|y)K.$$
\end{lem}
\begin{proof}
Clearly, $[x,y]$ has zero weight so $[x,y]\in\fhh$.
The invariance of $(-|-)$ gives $([x,y]|z)=(x| [y,z])$
so $([x,y]|h)=0$ if $h\in\fhh\cap [\fhg,\fhg]$
and $([x,y]|D)=m (x|y)$. 
\end{proof}

\subsubsection{}\label{cHN}
One has decompositions 
$$\fn_-=\cN^-_-\oplus\cH_-\oplus\cN^+_-,\ 
\fn=\cN^-\oplus\cH\oplus\cN^+,$$
where
$$\cH_-:=\fn_-\cap\tilde{\cH},
\ \ \cN^{\pm}_-:=\fn_-\cap\tilde{\cN}^{\pm},\ \ \cH:=\fn_-\cap\tilde{\cH},
\ \ \cN^{\pm}:=\fn\cap\tilde{\cN}^{\pm}.$$

If $\fhg$ is not of the type $A(2k,2l)^{(4)}$, $\cH_-$ is a (even) 
commutative countably dimensional Lie algebra.

\subsection{Verma modules}
Set $\fhb:=\fhh+\fhn$.
For $\lambda\in\fhh^*$ let $M(\lambda)$ be
the Verma module of the highest weight $\lambda$ and
let $v_{\lambda}$ be the canonical generator of $M(\lambda)$.
A Verma module $M(\lambda)$ has a unique maximal proper submodule
$M'(\lambda)$; denote by $L(\lambda)$ the simple quotient
 $L(\lambda):=M(\lambda)/M'(\lambda)$. 

\subsubsection{}\label{subV}
The invariant form $(-|-)$  gives rise to a quadratic Casimir
element which acts on $M(\lambda)$ by the scalar 
$(\lambda,\lambda+2\hrho)$. As a result, if
$L(\lambda-\nu)$ is a subquotient of $M(\lambda)$
then $2(\lambda+\hrho,\nu)=(\nu,\nu)$.

\subsubsection{}\label{mult:}
Verma modules do not admit Jordan-H\"older series, since
some Verma modules have an infinite length. 
However, so-called {\em local series}
introduced in~\cite{DGK} are  nice 
substitution for Jordan-H\"older ones. A series of weight modules
$N=N_0\supset N_1\supset\ldots\supset N_m=0$ is called 
{\em local at $\nu\in\fhh^*$}
if either $N_i/N_{i+1}\cong L(\lambda_i)$ for some $\lambda_i\geq \nu$ 
or $(N_i/N_{i+1})_{\mu}=0$ for all $\mu\geq\nu$.
This allows to define the {\em multiplicity} $[N:L(\lambda)]$
as the number of $i$ such that $N_i/N_{i+1}\cong L(\lambda)$
for a series local at some $\nu\leq \lambda$.

\subsection{Projections $\HC$ and $\HC_{\pm}$}\label{HC-}
We identify $\Uhh$ with $\cS(\fhh)$.
The triangular decomposition $\fhg=\fhn_-\oplus\fhh\oplus\fhn$ 
gives rise to the Harish-Chandra projection 
$\HC:\cU(\fhg)\to\cU(\fhh)=\cS(\fhh)$
along the kernel $\cU(\fhg)\fhn^++\fhn^-\cU(\fhg)$.
The restriction of $\HC$ to $\cU(\fhg)^{\fhh}$ 
is an algebra homomorphism.

\subsubsection{Algebra $\cS$}\label{cS}
Set
$$\cS:=\cU(\cH_-).$$
The algebra $\cS$ inherits the grading
$\cS=\oplus_{s\geq 0} \cS_{-s\delta}$.

Let $\fhg\not=A(2k,2l)^{(4)}$. Then
$\cH_-$ is pure even and commutative  so 
$\cS$ is the symmetric algebra of $\cH_-$.
In this case  $\cS$ has a natural grading
of symmetric algebra $\cS=\oplus \cS^k$,
where $\cS^k$ is spanned by the monomials 
of degree $k$. We set
$\cS^{\geq k}:=\sum_{j=k}^{\infty} \cS^j$.

\subsubsection{Projection $\HC_{\pm}$}\label{HC--}
The  decomposition $\fhn_-=\cN^+_-\oplus\cH_-\oplus\cN^-_-$
gives rise to the projections $\HC_{\pm}:\cU(\fhn_-)\to \cS$,
where
$$\Ker\HC_+=\cU(\fhn_-)\cN^-_-+\cN^+_-\cU(\fhn_-),\ \ \
\Ker\HC_-=\cU(\fhn_-)\cN^+_-+\cN^-_-\cU(\fhn_-).$$
Recall that $\fh'=\fhh\cap[\fhg,\fhg]$ and $\cH_-=\fhn_-^{\fh'}$.
Observe that $(\Ker\HC_+)^{\fh'}=\cU(\fhn_-)\cN^-_-=\cN^+_-\cU(\fhn_-)$
and similarly $(\Ker\HC_+)^{\fh'}=\cU(\fhn_-)\cN^+_-=\cN^-_-\cU(\fhn_-)$.
Thus the restriction of  $\HC_{\pm}$ to  $\cU(\fhn_-)^{\fh'}$ 
($\fh'=\fhh\cap[\fhg,\fhg]$)
are algebra homomorphism.

A Verma module $M(\lambda)$ can be canonically identified
with $\cU(\fhn_-)$. We denote by
$\HC_{\pm}: M(\lambda)\to\cS$ the corresponding
maps obtained via this identification.
Set
$$\breve{M}(\lambda):=\cU(\fhn_-)^{\fh'}v_{\lambda}$$
and note that $\breve{M}(\lambda)=\sum_s M(\lambda)_{\lambda-s\delta}$.
For the maps $\HC_{\pm}: M(\lambda)\to\cS$ the above observation gives
\begin{equation}\label{kerhcl}
(\Ker\HC_+\cap \breve{M}(\lambda))\subset
\cN^+_- M(\lambda),\ \ \ \ 
(\Ker\HC_-\cap \breve{M}(\lambda))\subset
\cN^-_- M(\lambda).
\end{equation}

\subsection{Shapovalov form}
Let us recall the standard construction of the Shapovalov form for $\fhg$.

\subsubsection{Choice of antiautomorphism.}
Call a linear endomorphism $\sigma$
of a superalgebra a ``naive'' antiautomorphism
if $\sigma$ is invertible and $\sigma([xy])=[\sigma(y),\sigma(x)]$.
By construction, any Kac-Moody  superalgebra admits
a ``naive'' anti-involution $\sigma$ which preserves the elements
of a Cartan subalgebra ($\sigma$ acts  on Cartan generators by
$\sigma(e_i)=f_i$).  One has
$$\sigma^2=\id,\  \sigma|_{\fhh}=\id,\ \sigma(\fhn)=\fhn_-,\ 
(\sigma(b)|\sigma(a))=(a|b).$$
All the above properties except the last one follows from
the definition of $\sigma$. To verify the last formula
note that $(a|b)'':=(\sigma(b)|\sigma(a))$ is an even invariant
bilinear form and both $(-|-), (-|-)''$ induce maps $\varphi,\varphi'':\fhg\to
\fhg^*$. The kernel of $\varphi-\varphi''$ is an ideal containing
$\fhh$ since $\sigma|_{\fhh}=\id$. Hence $\varphi=\varphi''$
and thus $(\sigma(b)|\sigma(a))=(a|b)$ as required.

\subsubsection{}\label{cosi}
Define a form $S: \cU(\fhn_-)\otimes \cU(\fhn_-)\to\cS(\fhh)$ by setting
$S(x,y):=\HC(\sigma(x)y)$. Using the natural
identification of a Verma module $M(\lambda)$  
with $\cU(\fhn_-)$, one easily sees that the maximal proper 
submodule $M'(\lambda)$ of $M(\lambda)$ coincides
with the kernel of the evaluated form 
$S(\lambda):\cU(\fhn_-)\otimes \cU(\fhn_-)\to \mathbb{C}$.

Notice that $S(x,y)=0$ if $x$ and $y$ have distinct weights. Thus
$S=\sum_{\nu\in \hat{Q}^+} S_{\nu}$ where $S_{\nu}$ is the restriction
of $S$ to $\cU(\fn_-)_{-\nu}\otimes \cU(\fn_-)_{-\nu}$.
By the above, $\dim L(\lambda)_{\lambda-\nu}=\codim\ker S_{\nu}(\lambda)$.

\subsubsection{}\label{det}
The determinant 
of $S_{\nu}$ is defined up to a non-zero scalar factor and
is computed in~\cite{kk}, \cite{gk}:
$$\det S_{\nu}(\lambda)=\prod_{m=1}^{\infty}
\prod_{\gamma\in \hat{\Delta}^+}
\phi_{m\gamma}^{(-1)^{(m-1)p(\gamma)}\tau(\nu-m\gamma)}(\lambda),$$
where $\hat{\Delta}^+$ stands for the multiset of positive roots, 
$\phi_{m\gamma}(\lambda)=2(\lambda+\hat{\rho},m\gamma)-(m\gamma,m\gamma)$,
and $\tau:\ \hat{Q}^+\to \mathbb{Z}_{\geq 0}$ stands for
the Kostant partition function given by 
$\tau(\nu)=\dim\cU(\fhn)_{\nu}=\dim M(\lambda)_{\lambda-\nu}$.

Note that for any $\gamma\in\hat{\Delta}_{im}, m\geq 1$ 
one has $\phi_{m\gamma}(\lambda)=0$ iff $(\lambda+\hat{\rho},\delta)=0$.

\subsubsection{}\begin{rem}{remSorth}
Note that $\sigma(\cN^+)=\cN^-_-,\ \sigma(\cN^-)=\cN^+_-$ because
$\sigma|_{\fhh}=\id$. From~(\ref{kerhcl}) it follows
that $\breve{M}(\lambda)^{\cN^{\pm}}$ is orthogonal to
$\Ker\HC_{\mp}$ with respect to the Shapovalov form:
\begin{equation}\label{orth}
S(\lambda)\bigl(\breve{M}(\lambda)^{\cN^+},\Ker\HC_- \bigr)=
S(\lambda)\bigl(\breve{M}(\lambda)^{\cN^-},\Ker\HC_+ \bigr)=0.
\end{equation}
\end{rem}

\subsection{Characters}\label{Omega}
For a diagonalizable $\fhh$-module $M$ we denote by $\Omega(M)$
the set of weights of $M$ and by $M_{\mu}$ 
the weight space of weight $\mu$. 
We say that a module $M$ admits a {\em character} 
$\ch M=\sum_{\mu}\dim M_{\mu} e^{\mu}$
if $M$ is a diagonalizable $\fhh$-module with 
the finite dimensional weight spaces.

If $v\in M$ is a weight vector
we denote its weight by $\wt v$;
for $u\in\Ug$ we denote by $\wt u$ its weight with respect to 
the adjoint action of $\fhh$.

\subsubsection{}\label{xy}
For a finite set $X:=\{\lambda_i\}_{i=1}^r\subset \fhh$ let
$C_X$
be the collection of elements of the form 
$\sum_{i=1}^r \sum_{\mu<\lambda_i} c_{\mu} e^{\mu}$,
where $c_{\mu}\in \mathbb{Z}, c_{\mu}\geq 0$. Let $C$
be the union of $C_X$. 
Note that $x,y\in C$  implies $x+y,\ xy\in C$ and 
introduce the partial order on $C$ by $x\geq y$ if $x-y\in C$.
In all our examples $\ch M$ belongs to $C$.

\subsubsection{}
\label{mu}
Set $M_{>\mu}:=\sum_{\nu>\mu} M_{\nu}$ and define $M_{\geq\mu}$ similarly.
For any map $f:M\to N$ we denote by $f_{\mu}$ (resp., $f_{>\mu},
f_{\geq\mu}$) the restriction
of $f$ to $M_{\mu}$ (resp., $M_{>\mu}, M_{\geq\mu}$).

\subsubsection{}\label{chrzeta}
Recall that $\fhg=[\fhg,\fhg]\oplus \mathbb{C}D$.
As $[\fhg,\fhg]$-modules for any $s$
one has $M(\lambda)\cong M(\lambda-s\delta)$ 
and so $M'(\lambda)\cong M'(\lambda-s\delta)$. As a consequence,
$\ch L(\lambda-s\delta)=e^{-s\delta}\ch L(\lambda)$.

\subsection{Generic critical weights}\label{defgn}
Set
$$\Lambda:=\{\lambda\in\fhh|\ [M(\lambda):L(\lambda-\nu)]=0\ \text{ for }
\nu\not\in\mathbb{N}\delta\},\ \ \
\Lambda_{crit}=\{\lambda\in\Lambda|\ (\lambda+\hat{\rho},\delta)=0\}.$$
From the formula for Shapovalov determinants (see~\ref{det}) we
see that $\lambda\in\Lambda\setminus\Lambda_{crit}$ iff 
$M(\lambda)$ is simple. For $\lambda\in\Lambda_{crit}$
the Verma module $M(\lambda)$ is not irreducible and all
its subquotients are of the form $L(\lambda-s\delta)$ for $s\geq 0$.

By~\cite{kk}, Thm. 2 
$\lambda\in\Lambda_{crit}$ iff 
$(\lambda+\hrho,\delta)=0$ and $\lambda$ is not a root of other
factors of Shapovalov determinants, i.e
$\phi_{m\gamma}(\lambda)\not=0$ for $\gamma\not\in\Delta_{im}$.

Call $\lambda\in\Lambda_{crit}$ a {\em generic critical weight} if
$\lambda\in\Lambda_{crit}$.

\subsubsection{}
We will use the following lemma.

\begin{lem}{primcrit}
Take $\lambda\in\Lambda$. If
$v'\in M(\lambda)_{\lambda-s\delta}$ is such that
$(\cN^++\cH)v'=0$ or $(\cN^-+\cH)v'=0$
then $v'$ is singular.
\end{lem}
\begin{proof}
Assume that $(\cN^++\cH)v'=0$.
Since $\fhn=(\cN^++\cH)\oplus {\cN^-}$ one has
$\cU(\fhn)v'=\cU({\cN^-})v'$. 
The weight vectors in $\cU({\cN^-})v'$
which are not proportional to $v'$ have weights
of the form $\lambda-\mu$ where $\mu\not\in\mathbb{N}\delta$;
these vectors are not singular, since $\lambda\in\Lambda$.
Since $\cU(\fhn)v'$ contains at least one singular weight vector,
$v'$ is singular.

The proof for $(\cN^-+\cH)v'=0$ is similar.
\end{proof}

\subsubsection{}
\begin{rem}{remgn}
Recall that each factor of a Shapovalov determinant corresponds
to a hyperplane $\gamma_{k\alpha}$ in $\fhh^*$;
 the set of generic critical weights
consists of the points of the hyperplane $(\lambda+\hrho,\delta)=0$
which do not lie on other hyperplanes $\gamma_{k\alpha}$. From the formula 
given in~\ref{det} one sees that
$\lambda$ is a generic critical weight iff
$$(\lambda+\hrho,\delta)=0\ \& \ (\lambda+\hrho,\alpha)\not
=k(\alpha,\alpha),
\text{ where } \alpha\in\hat{\Delta}^+,\ 
k\in\mathbb{Z}_{>0}\ \& \
k\text{ is odd if }\alpha \text{ is odd}.$$
\end{rem}

\section{Towards the proof of Theorems~\ref{thm01},\ref{thm02}}
\label{sectorm}
In this section we reduce Theorems~\ref{thm01},\ref{thm02}
to the assertions~\ref{reform} (A), (B).

\subsection{Algebra structure on $\breve{M}(\lambda)^{\fhn}$}
\label{algstr}
Recall that $\breve{M}(\lambda)=\sum_s M(\lambda-s\delta)$.
Identify $\End_{[\fhg,\fhg]}(M(\lambda))$ with 
$\breve{M}(\lambda)^{\fhn}$ via the map
 $\phi\mapsto \phi(v_{\lambda})$ and
endow  $\breve{M}(\lambda)^{\fhn}$ with the algebra structure via
 this bijection.

\subsubsection{}
\begin{lem}{lemrst}
The restrictions $\HC_{\pm}: \breve{M}(\lambda)^{\fhn}\to\cS$ are algebra
 homomorphisms.
\end{lem}
\begin{proof}
 Let $\iota$ stands for the natural identification
$M(\lambda)$ with $\cU(\fhn_-)$. Clearly,
the algebra structure on $\breve{M}(\lambda)^{\fhn}$
  is compatible with $\iota$, i.e. $\iota(x)\iota(y)=\iota(xy)$.

Notice that $\iota$ maps $\breve{M}(\lambda)$
to $\cU(\fhn_-)^{\fh'}$, where
$\fh'=\fhh\cap[\fhg,\fhg]$. Now the assertion follows
from the fact that the restrictions of $\HC_{\pm}$ to
$\cU(\fhn_-)^{\fh'}$ is an algebra homomorphism.
\end{proof}

\subsection{}
\begin{claim}{reform}
Take $\fhg\not=A(2k,2l)^{(4)}$ and $\lambda\in\Lambda_{crit}$.
The following assertions imply Theorems~\ref{thm01}, \ref{thm02}
and the formula~(\ref{cor01}):
 
(A) The restrictions of
$\HC_{\pm}$ to the set of singular vectors of $M(\lambda)^{k}$ 
contain $\cS^{\geq k}$.

(B) $\forall k,s\ \ [M(\lambda)^k:L(\lambda-s\delta)]
=\dim\cS^{\geq k}_{-s\delta}$.
\end{claim}

We prove this claim in~\ref{pfi}---\ref{phiv} below.

\subsubsection{} \label{pfi}
\Thm{thm01} (i) claims that $M(\lambda)$ 
does not have subsingular vectors, i.e.
\begin{equation}\label{MM}
\forall\mu\in\fhh^*\ \ \ \ \ \ [M(\lambda):L(\mu)]=
\dim M(\lambda)^{\fhn}_{\mu}\end{equation}
and that any submodule of 
$M(\lambda)$ is generated by its singular vectors.
Combining (A) and (B) for $k=0$ we obtain~(\ref{MM}).

Let us show that~(\ref{MM}) implies that
any submodule of $M(\lambda)$ is generated by its singular vectors.
Let $\cY$ be the set of modules in category $\cO$ which do
 not have subsingular vectors:
$\cY:=\{M\in\cO|\ \forall\mu\ 
[M: L(\mu)]=\dim M^{\fhn}_{\mu}\}$.
Take $M\in\cY$ and let $N$ is a submodule of $M$. For any $\mu\in\fhh^*$
one has
$$[M: L(\mu)]=[N: L(\mu)]+[M/N: L(\mu)]\geq 
\dim N^{\fhn}_{\mu}+\dim (M/N)^{\fhn}_{\mu}\geq \dim M^{\fhn}_{\mu}.$$
Therefore $N,\ M/N\in\cY$ and, moreover, 
$\dim M^{\fhn}_{\mu}=\dim N^{\fhn}_{\mu}+\dim (M/N)^{\fhn}_{\mu}$.

By~(\ref{MM}), $M(\lambda)\in\cY$.
Let $M$ be any submodule of $M(\lambda)$ and
let $N$ be the submodule generated by $M^{\fhn}$.
By above $M\in\cY$ and
$\dim (M/N)^{\fhn}_{\mu}=0$ for any $\mu$. Hence $M/N=0$ as required.

\subsubsection{}\label{phi2}
\Thm{thm01} (ii) asserts that $\HC_{\pm}$
induce bijections
$\HC_{\pm}: M(\lambda)^{\fhn}\iso\cS$
and that these bijections are algebra isomorphisms.
By (A) the maps $\HC_{\pm}: M(\lambda)^{\fhn}\to\cS$
are surjective. Combining~(\ref{MM}) with (B) for $k=0$
we get $e^{-\lambda}\ch M(\lambda)^{\fhn}=\ch \cS$
and this implies the injectivity.
The fact that these bijections are algebra isomorphisms
follows from~\Lem{lemrst}.

\subsubsection{} 
\Thm{thm01} (iii) asserts that the submodule generated by
a singular vector $v\in M(\lambda)$ is isomorphic
a Verma module.

Let $\phi\in \End_{[\fhg,\fhg]}(M(\lambda))$ 
be such that $\phi(v_{\lambda})=v$. If $\Ker\phi$ is non-zero,
it contains a singular vector $w$.
Then $0=\HC_+(\phi(w))=\HC_+(w)\HC_+(v)$ by~\Lem{lemrst}.
By~\ref{phi2}, $\HC_+(w), \HC_+(v)$ are non-zero. Since $\cS$
is a polynomial algebra this gives a contradiction.
Hence $\Ker\phi=0$ and this implies (iii).

\subsubsection{}\label{pfi3}
Let us deduce the formula~(\ref{cor01}):
$\ch N=\ch L(\lambda)\cdot \ch \HC_+(N^{\fhn}),$
where $N$ is a submodule of $M(\lambda)$. 

By~\ref{pfi}, $[N:L(\mu)]=\dim N^{\fhn}_{\mu}$.
Since the weights of singular vector are of the form
$\lambda-s\delta$, one has
$\ch N=\ch L(\lambda)\cdot (e^{-\lambda}\ch N^{\fhn})$.
By~\ref{phi2},   $e^{-\lambda}\ch N^{\fhn}=\ch \HC_+(N^{\fhn})$
and this gives the required formula.

\subsubsection{}\label{phiv}
Recall that~\Thm{thm02}  means that $\HC_{\pm}$
induce bijections
$\HC_{\pm}: \bigl(M(\lambda)^k\bigr)^{\fhn}\iso\cS^{\geq k}$.
By~\ref{phi2}, the restrictions of
$\HC_{\pm}$ to $\bigl(M(\lambda)^k\bigr)^{\fhn}$ are injective;
by (A) the images contain $\cS^{\geq k}$.
Finally, combining~(\ref{MM}) with (B) we get 
$e^{-\lambda}\ch\bigl(M(\lambda)^k\bigr)^{\fhn}=\ch\cS^{\geq k}$.
Hence the images coincide with $\cS^{\geq k}$ as required.
\qed

\subsection{}
\begin{cor}{corH} For $\lambda$ at the critical level one has
$\ch  L(\lambda)\ch \cS\leq \ch M(\lambda)$.
\end{cor}
\begin{proof}
By (A) for $\lambda\in\Lambda_{crit}$ one has
$e^{-\lambda}\ch M(\lambda)^{\fhn}\geq \ch\cS$.
Since $\Lambda_{crit}$ is Zariski dense at the critical level, 
$\dim M(\lambda)^{\fhn}_{\lambda-s\delta}\geq\dim\cS_{s\delta}$
for any $\lambda$ at the critical level. In the
notation of~\ref{mult:} we have
$$\begin{array}{rl}
\ch M(\lambda)&=\sum_{\nu\in\hat{Q}^+}
[M(\lambda):L(\lambda-\nu)]\ch L(\lambda-\nu)\geq 
\displaystyle\sum_{s=0}^{\infty} [M(\lambda):L(\lambda-s\delta)]
\ch L(\lambda-s\delta)\\
&=^{\text{by~\ref{chrzeta}}}
\displaystyle\sum_{s=0}^{\infty} [M(\lambda):L(\lambda-s\delta)]
e^{-s\delta}\ch L(\lambda)\geq 
\ch L(\lambda)\displaystyle\sum_{s=0}^{\infty}\dim\cS_{-s\delta}
e^{-s\delta}\\
&=\ch  L(\lambda)\ch \cS. \qed
\end{array} $$
\end{proof}

\section{Jantzen filtration}\label{sectjan}
In this section we recall the construction
of the Jantzen filtration (see~\cite{Jan}) and then describe
this filtration for a vacuum module over a Heisenberg algebra.

\subsection{Notation}\label{notloc}
Let $A$ be the localization of $\mathbb{C}[x]$ by the ideal
generated by $x$: $A:=\mathbb{C}[x]_{(x)}$.
For a Lie algebra $\fp$ and a $\fp$-module $N$ set
$\fp_A:=\fp\otimes_{\mathbb{C}} A$, $N_A:=N\otimes_{\mathbb{C}} A$.

\subsubsection{Choice of $\xi$}
Choose an element
$\xi\in\fhh^*$ such that for any $\nu\in\hat{Q}^+$
the Shapovalov determinant $\det S_{\nu}(\lambda+t\xi)$,
viewed as a polynomial in $t$, is not identically equal to zero.

\subsubsection{}\label{MuuR}
Introduce a Verma module $M(\lambda+x\xi)$ over $\fhg_A$ as follows.

Define the action of  $(\fn+\fh)$ on  $A$: 
 $\fn$ acts trivially and
$h\in\fh$ acts by the multiplication to 
$(\lambda+x\xi)(h)=\lambda(h)+x\xi(h)$.
Now $M(\lambda+x\xi)$ is the following $\fhg_A$ module
$$M(\lambda+x\xi):=
\Ind_{\fn+\fh}^{\fg} A.$$

 Set
$$\breve{M}(\lambda+x\xi)=\sum_{s=0}^{\infty}
M(\lambda+x\xi)_{\lambda+x\xi-s\delta}.$$

\subsection{Definition of Jantzen filtration}\label{shapR}
The Shapovalov form $S:\cU(\fhn_-)\otimes\cU(\fhn_-)\to\cS(\fhh)$
induces the Shapovalov form 
$S_A: \cU(\fhn_{-,A})\otimes\cU(\fhn_{-,A})\to\cS(\fhh_A)$
and  the evaluated Shapovalov form 
$S_A(\lambda+x\xi): M(\lambda+x\xi)\otimes M(\lambda+x\xi)\to A$.
For  $k\in\mathbb{Z}_{\geq 0}$, set
$$\cF(M(\lambda+x\xi))^k:=\{v\in M(\lambda+x\xi)|\ 
   S_A(\lambda+x\xi)(v,v')\in t^kA \ \ \forall v'\}.$$
This defines a decreasing filtration. It is easy to see 
that each term $\cF(M(\lambda+x\xi))^k$ is a submodule
of $M(\lambda+x\xi)$. 
The weight spaces of $M(\lambda+x\xi)$ are free of finite rank
$A$-modules and the determinant of the restriction
of $S_A(\lambda+x\xi)$ to the weight space $\lambda+x\xi-\nu$ i
equal (up to a non-zero factor in $A$) to
$\det S_{\nu}(\lambda+x\xi)$
and this is non-zero due to the condition on $\xi$. As a result, 
$\cap_{k=0}^{\infty} \cF^k(M(\lambda+x\xi))=0$.

Specializing this filtration at $x=0$ we obtain
the Jantzen filtration $M^k(\lambda)$ on $M(\lambda)$.
The terms $M^k(\lambda)$
form a decreasing filtration by submodules of $M(\lambda)$
having zero intersection. One readily sees that 
$M^0(\lambda)=M(\lambda)$ and that
$M^1(\lambda)$ coincides with the maximal proper
submodule of $M(\lambda)$.

Define $\upsilon_x: A\to \mathbb{N}$ by
$a\in (x)^{\upsilon_x(a)},
\ a\not\in (x)^{\upsilon_x(a)+1}$ for any $a\in A$.
The Jantzen sum formula (\cite{Jan}, Lem.3) claims that
$$\forall\mu\ \ \ 
\sum_{k\geq 1}\dim M(\lambda)^k_{\lambda-\mu}=
\upsilon_x\bigl(\det S_{A;\mu}(\lambda+x\xi)\bigr).$$
Observe that for each $\mu$ one has $M(\lambda)^k_{\mu}=0$ for $k>>0$.

\subsubsection{}
Let $\lambda$ be a generic critical weight.
For $\fhg\not=A(22k,2l)^{(4)}$ the set $\hat{\Delta}^+_{im,\ol{1}}$
is empty and the Jantzen sum formula gives
\begin{equation}
\label{sumFA'}
\sum_{k\geq 1}\ch M(\lambda)^k=\ch M(\lambda)\cdot
\sum_{m=1}^{\infty}\sum_{\alpha\in\hat{\Delta}^+_{im}}e^{-m\alpha}=
\ch M(\lambda)\cdot
\sum_{\alpha\in\hat{\Delta}^+_{im}}\frac{e^{-\alpha}}{1-e^{-\alpha}}.
\end{equation}

\subsection{Jantzen filtration for modules over a Heisenberg algebra}
\label{vac}
Let $\fhg\not=A(2k,2l)^{(4)}$. In this case 
$\cH_-+\mathbb{C}K+\cH$ is a Heisenberg algebra. 
For $k\in\mathbb{C}$ let $V^k$ be
 a vacuum module over
$\cH_-+\mathbb{C}K+\cH$: 
$$V^k=\cU(\cH_-+\mathbb{C}K+\cH)\otimes_{\cU(\mathbb{C}K+\cH)} 
\mathbb{C}v_k,\ \ Kv_k=kv_k, \cH v_k=0.$$
View $V^k$ as a graded module via the natural identification
of $V^k$ with $\cS=\cU(\cH_-)$; note that
 the homogeneous components of $V^k$
are finite dimensional.

Define the Shapovalov form $\cS\otimes\cS\to \mathbb{C}[K]$
as in~\ref{cosi}.  It is easy to see that 
$V^k$ is simple for $k\not=0$ and that for $V^0$ the Jantzen filtration 
$\cF^k(V^0)$ identifies with the adic filtration 
$\{\cS^{\geq k}\}$ of $\cS=\cU(\cH_-)$. Set 
$\ch_q\cS=\prod_{\alpha\in\hat{\Delta}^+_{im}}(1-qe^{-\alpha})^{-1}$.
Then
$$\sum_{k\geq 1}\ch \cS^{\geq k}=
\frac{\partial \ch_q\cS}{\partial q}|_{q=1}.$$

Thus 
\begin{equation}
\label{sumFS}
\sum_{k\geq 1}\ch \cF^k(V^0)=\sum_{k\geq 1}\ch \cS^{\geq k}=\ch \cS
\sum_{\alpha\in\hat{\Delta}^+_{im}}\frac{e^{-\alpha}}{1-e^{-\alpha}}.
\end{equation}

\subsubsection{}
Let $\lambda$ be a generic critical weight.
Combining~(\ref{sumFA'}) and~(\ref{sumFS})
 we get
\begin{equation}
\label{sumMS}
\sum_{k\geq 1}\ch M(\lambda)^k =\frac{\ch M(\lambda)}{\ch \cS}
\sum_{k\geq 1}\ch \cS^{\geq k}.\end{equation}

Recall that the 
Kac-Kazhdan character formula states that $\ch L(\lambda)=
\frac{\ch M(\lambda)}{\ch \cS}$ 
and that~\Thm{thm02}  states that 
$\ch M(\lambda)^k=\ch L(\lambda)\ch \cS^{\geq k}$ for all $k$.
We will use the formula~(\ref{sumMS}) in the proof of~\ref{reform}
 (B).

\section{A technical lemma}\label{secttech}
In this section $\fhg$ is an affine Lie superalgebra with a symmetrizable
Cartan matrix (we do not exclude the case $A(2k,2l)^{(4)}$).
We prove~\Lem{lemVhm}; this lemma is
used in the explicit construction of a singular vector in
Sect.~\ref{sectexp}.

In this section $\fq,\ \fq_-$ are any subalgebras of $\fhn$ satisfying
$$\cN^+\subseteq \fq\subseteq \cN^++\cH,\ \ 
 \ \cN^-\subseteq \fq^-\subseteq \cN^-+\cH.$$

\subsection{Notation}\label{notlam}
Let $h'\in\fhh$ be the defining element for the triangular
decomposition $[\fhg,\fhg]=\tilde{\cN^-}\oplus
\tilde{\cH}\oplus\tilde{\cN^+}$ (see~\ref{trfhg}). 
Recall that 
 $\Omega(\tilde{\cN^+})=\{\alpha\in\hat{\Delta}|\ \Ree
\alpha(h')>0\}$.
Set $\hat{Q}^-:=\Omega(\cU(\fhn^-))=-\hat{Q}^+$ and
$$\hat{Q}^-_+:=\{\nu\in \hat{Q}^-|\ \Ree\nu(h')>0\},\ \ 
\hat{Q}^-_-:=\{\nu\in \hat{Q}^-|\ \Ree\nu(h')<0\}.$$

Retain notation of~\ref{defgn}.

\subsection{}
\begin{lem}{chmlem2}
For any $\lambda\in\Lambda$ one has
$$\begin{array}{l}
H^m(\fq, L(\lambda))_{\lambda+\mu}=0\ \text{ and }\ 
H^m(\fq, M(\lambda))_{\lambda+\mu}=0
\ \text{ for } m=0,1\ \text{ and } \mu\in \hat{Q}^-_-,\\
H^m(\fq^-, L(\lambda))_{\lambda+\mu}=0\ \text{ and }\ 
H^m(\fq^-, M(\lambda))_{\lambda+\mu}=0
\ \text{ for } m=0,1\ \text{ and } \mu\in \hat{Q}^-_+.\end{array}$$
\end{lem}
\begin{proof}
Recall that $\sigma(\cN^+)=\cN^-_-$.
 \Lem{chmlem1} reduces the first
formula $H^m(\fq, L(\lambda))_{\lambda+\mu}=0$
to the inclusions
$$L(\lambda)_{\lambda+\mu}\subset \cN^-_-L(\lambda),\ \ \ 
M'(\lambda)_{\lambda+\mu}\subset \cN^-_- M'(\lambda),$$
where $M'(\lambda)$ is the maximal proper submodule of $M(\lambda)$.

Recall that $\Ree((s\delta+\mu)h')<0$.
Observe that if $\nu$ is a weight of $\cU(\cN^+_-+\cH_-)$ then
$\Ree\nu(h')\geq 0$; thus
$\nu$ is not of the form
$s\delta+\mu$ for $s\in\mathbb{Z}$:
$$\Omega(\cU(\cN^+_-+\cH_-))\cap 
\bigl(\mu+\mathbb{Z}\delta\bigr)
=\emptyset.$$
Both inclusions can be easily deduced from this observation.

Indeed, writing $\fhn_-=\cN^-_-\oplus(\cN^+_-+\cH_-)$ we get
$\cU(\fhn_-)=\cU(\cN^+_-+\cH_-)\oplus\cN^-_-\cU(\fhn_-)$.
One has $L(\lambda)=\cU(\fhn_-)\ol{v}_{\lambda}$, where $\ol{v}_{\lambda}$
is the highest weight vector in $L(\lambda)$.
By the above observation, $L(\lambda)_{\lambda+\mu}$ does not meet 
$\cU(\cN^+_-+\cH_-)\ol{v}_{\lambda}$ and this gives the first inclusion.
For the second inclusion, recall that
 $M'(\lambda)=0$ if the level of $\lambda$ is not critical
and that $M'(\lambda)=\cU(\fhn_-)N$,
where $N:=\sum_{j\geq 1}M'(\lambda)_{\lambda-j\delta}$ if
$\lambda$ has the critical level.  By the above observation,
the space $\cU(\cN^+_-+\cH_-)N$ does not meet $M(\lambda)_{\lambda+\mu}$.
Hence $M'(\lambda)_{\lambda+\mu}\subset \cN^-_-M'$ and
the second inclusion follows.

Recall that
$M(\lambda)$ has a local series at $\lambda+\mu$ with simple quotients
$L(\lambda-s_i\delta)$, where $s_i$ are non-negative integers. 
Applying the first formula
to the pair $\lambda':=\lambda-s_i\delta,\ \mu':=\mu+s_i\delta$ we
get $H^m(\fq, L(\lambda-s_i\delta))_{\lambda+\mu}=0$ for $m=0,1$ and
any index $i$. Using
the long exact sequence of Lie algebra cohomology we obtain
$H^m(\fq, M(\lambda))_{\lambda+\mu}=0$.

The proof for $\fq^-$ are  similar.
\end{proof}

\subsection{}
Retain notation of~\ref{notloc}.

\begin{cor}{chmloc}
For $\lambda\in\Lambda$ one has
$$\begin{array}{l}
H^m(\fq_A, M(\lambda+x\xi))_{\lambda+x\xi+\mu}=0
\ \text{ for } m=0,1\ \text{ and } \mu\in \hat{Q}^-_-,\\
H^m(\fq^-_A, M(\lambda+x\xi))_{\lambda+x\xi+\mu}=0
\ \text{ for } m=0,1\ \text{ and } \mu\in \hat{Q}^-_+.\end{array}$$
\end{cor}
\begin{proof}
Set $R:=\mathbb{C}[x]$ and let $M$
be a Verma $\fhg_R$-module of highest weight $\lambda+x\xi$.
Clearly,  $M(\lambda+x\xi)=M\otimes_R A$.
For $c\in\mathbb{C}$ denote by $M(c)$ the evaluation of
$M$ at the point $t=c$. Clearly, $M(c)$ is a Verma $\fhg$-module 
isomorphic to $M(\lambda+c\xi)$. The choice of $\xi$ ensures
that $(\lambda+c\xi)\in\Lambda$ for infinitely many values of $c$.
Combining Lemmas~\ref{chmlem2} and~\ref{lemcohoR}  we obtain the required
assertion.
\end{proof}

\subsection{}\label{Vhm}
Let $V_1,V_2$ be any diagonalizable $\fhh$-modules, i.e. 
$V_i=\sum_{\xi\in\fhh^*} V_{i,\mu}$. Call a linear map
$\varphi: V_1\to V_2$ a {\em weight homomorphism}
if there exists $\nu\in\fhh^*$ such that
$\varphi(V_{1,\mu})\subset V_{2,\mu+\nu}$ for all $\mu\in\fhh^*$.

\subsubsection{}
Let $W$ be an $(\fhh+\fq)$-module with finite dimensional weight
spaces (with respect to $\fhh$) such that for some $\mu\in\fhh^*$
one has $\dim W_{\mu}=1$ and
$$\Omega(W)\subset \{\mu+\nu|\ \nu\in\hat{Q}^-_-\cup\{0\}\}.$$
Let $w\in W_{\mu}$ be a non-zero vector.
Retain notation of~\ref{notloc}.

\subsubsection{}
\begin{lem}{lemVhm}
\begin{enumerate}
\item 
For each weight vector $w'\in \breve{M}(\lambda)$ 
satisfying $\fq w'=0$
there exists a unique weight $\fq$-homomorphism
$\psi: W\to M(\lambda)$ such that $\psi(w)=w'$.
\item
For each weight vector $w'\in \breve{M}(\lambda+x\xi)$ 
satisfying $\fq_A w'=0$
there exists a unique weight $\fq_A$-homomorphism
$\psi: W_A\to M(\lambda+x\xi)$ such that $\psi(w)=w'$.
\end{enumerate}
\end{lem}
\begin{proof}
Let us prove (ii); the proof of (i) is the same.
Recall that  $\Omega(\fq)\subset \hat{Q}^+$ and so the action of $\fq$
raises the weight. As a result,
$W$  admits an increasing $(\fhh+\fq)$-filtration
$\{W^k\}_{k\geq 0}$ 
with one-dimensional factors and $W^0=W_{\mu}$. Then
$\{W^k_A\}_{k\geq 0}$ is an increasing $(\fhh+\fq)_A$-filtration
with factors of rank one over $A$.

Shift the $\fhh_A$-action on $M(\lambda+x\xi)$ 
by $\wt w-\mu$, i.e. set
$h.v=(h-(\wt w-\mu)(h))v$
for $h\in\fhh_A, v\in M(\lambda+x\xi)$.  
The shifted action is compatible with the action of
$\fq_A$; view  $M(\lambda+x\xi)$ as $(\fhh+\fq)_A$-module
with respect to this shifted action. Notice that
$\psi$ is an $(\fhh+\fq)_A$-homomorphism.
In the formulas below we use this shifted $(\fhh+\fq)_A$-module
structure on $M(\lambda+x\xi)$ and
$\Hom$ stands for  $\Hom_{(\fhh+\fq)_A}$.

Let $T^{\nu}$ be a one-dimensional $(\fhh+\fq)$-module of weight $\nu$
(i.e. $\fq T^{\nu}=0,\ h|_{T^{\nu}}=\nu(h)\id$).
The short exact sequence
$$0\to W^{k-1}_A\to W^{k}_A\to T^{\nu}_A\to 0$$ 
gives
$$\begin{array}{r}
0\to \Hom (T^{\nu}_A, 
M(\lambda+x\xi))=H^0(\fq_A,M(\lambda+x\xi))_{\nu'}
\to\Hom (W^k_A, M(\lambda+x\xi))\to\\
\to \Hom(W^{k-1}_A, M(\lambda+x\xi))
\to \Ext^1(T^{\nu}_A,M(\lambda+x\xi))=
H^1(\fq_A,M(\lambda+x\xi))_{\nu'},\end{array}$$
where $\nu':=\wt w+\nu=\lambda+x\xi-s\delta+\nu$ for some $s\geq 0$.
For $k>0$ one has $\nu\in \hat{Q}^-_-$ so either 
$\nu-s\delta\not\in\hat{Q}^-$
or $\nu-s\delta\in\hat{Q}^-_-$. \Cor{chmloc} gives
$$H^0(\fq_A,M(\lambda+x\xi))_{\nu'}=
H^1(\fq_A,M(\lambda+x\xi))_{\nu'}=0.$$
Therefore $\Hom (W^k_A, M(\lambda+x\xi))=
\Hom (W^{k-1}_A, M(\lambda+x\xi))$. As a result,
$\Hom (W^k_A, M(\lambda+x\xi))=\Hom (W^0_A, M(\lambda+x\xi))=
\Hom(Aw, M(\lambda+x\xi))$. The assertion follows.
\end{proof}

\subsubsection{}\label{Vhm-} 
The same arguments allow us to deduce 
the facts similar to~\Lem{lemVhm} for $\fq^-$, $\fq^-_A$:
$\fq$ should be substituted by $\fq^-$ and $\hat{Q}^-_-$
(in the formula for $\Omega(W)$)
by $\hat{Q}^-_+$.

\section{Module $M(\lambda+x\xi)$}\label{xxi}
In this section $\lambda\in\Lambda$.
We consider a Verma $\fhg_A$-module $M(\lambda+x\xi)$
introduced in~\ref{MuuR}.

\subsection{Content of the section}\label{uM}
We will use the notation $\HC_{\pm}$ for the maps
$\cU(\fhn_-)_A \to\cS_A$ as well as for
the corresponding maps
$M(\lambda+x\xi)\to\cS_A$ and their restrictions.
In the following proposition the maps
$\HC_{\pm}: \breve{M}(\lambda+x\xi)^{\cN^{\pm}}\to
\cS_{A}$ are the restrictions of 
$\HC_{\pm}: M(\lambda+x\xi)\to\cS_A$.

\subsubsection{}
\begin{prop}{propz1}
For $\lambda\in\Lambda$ the maps
$\HC_+: \breve{M}(\lambda+x\xi)^{\cN^+}\iso
\cS_{A},\ \ 
\HC_-: \breve{M}(\lambda+x\xi)^{\cN^-}\iso
\cS_{A}$ are bijective.
\end{prop}

\subsubsection{}\label{outz1}
The proposition is proven in~\ref{v+}, \ref{HC-1z}, \ref{injHC}
below. 
In~\ref{v+} we explicitly construct 
$v_{\pm}\in \breve{M}(\lambda+x\xi)^{\cN^{\pm}}$ 
satisfying $\HC_{\pm}(v_{\pm})=z$
for each $z\in\cH_-$. Then
in~\ref{HC-1z} we construct such vectors 
for any $z\in\cS$ and thus the surjectivity
of the maps $\HC_{\pm}:  \breve{M}(\lambda+x\xi)^{\cN^{\pm}}\to\cS_A$
is established.
Finally, in~\ref{injHC}, we show the injectivity of these maps.

\subsubsection{}
Eventhough we are interested in the case $\lambda\in\Lambda_{crit}$,
we consider $\lambda\in\Lambda$.
We use this twice: first,  the proof of injectivity
of $\HC_{\pm}:  \breve{M}(\lambda+x\xi)^{\cN^{\pm}}\to\cS_A$
is based on the surjectivity of 
$\HC_{\pm}:  \breve{M}(\lambda'+t\xi)^{\cN^{\pm}}\to\cS_A$,
where $t\in\mathbb{C}$ is such that $\lambda'=\lambda+t\xi$ lies in 
$\Lambda\setminus\Lambda_{crit}$. The second time we use 
this in~\Lem{lemququ}, where we prove that, for 
$\fhg\not=A(2k,2l)^{(4)}$, the images of vectors
$v_{\pm}$ in $M(\lambda)$
are singular if $\lambda\in\Lambda_{crit}$. 
The proof of this lemma is based on the results of~\ref{wvlambda}, where
we study the action of $\cH$ on
$v_{\pm}$ for $z\in\cH_-$.
We do not use the results of~\ref{wvlambda}
in this section and~\ref{wvlambda} is not a part of 
the proof of~\Prop{propz1}.

\subsection{}\label{outv+}
In~\ref{v+} we will construct 
$v_{\pm}\in \breve{M}(\lambda+x\xi)^{\cN^{\pm}}$ 
satisfying $\HC_{\pm}(v_{\pm})=z$
for $z\in\cH_-$.
Let  us describe the construction of 
$v_+$ for the non-twisted case (i.e. $\fhg=\fhg'$).
The constructions of $v_-$ is completely similar.
The twisted case differs by some technical details.

Take $\lambda\in\Lambda$ and $z=ht^{-m}$ ($h\in\fh,\ m>0$). Set 
$$\fq:=\fn+\sum_{s\geq 1} (\fh+\fn)t^s,\ \ 
N:=\mathbb{C}h+\sum_{s\geq 0} \fn t^s,\ \ 
V:=\mathbb{C}h+\sum_{0\leq s \leq m-1} \fn t^s.$$
Observe that $N$ is a  $\fq$-submodule of $\fhg$.
View $N^*$
as a $\fq$-module via the antiautomorphism $-\id$
and let $h^*\in N^*$ be the  element dual to  $h$: $h^*(h)=1,\ 
h^*(\sum_{s\geq 0} \fn t^s)=0$.
By~\Lem{lemVhm} there exists a unique $\fq_A$-homomorphism
$\psi: N^*_A\to M(\lambda+x\xi)$ such that $\psi(h^*)$
is the highest weight vector of $M(\lambda+x\xi)$.

For $a\in\fg$ set $T_m(at^s):=at^{s-m}$.
Let $\gamma: \fhg_A\otimes M(\lambda+x\xi)\to M(\lambda+x\xi)$ 
be the natural map $\gamma(u\otimes v)=uv$, and  let 
$\id'\in V_A\otimes N^*_A$ 
correspond to the identity map $\id_V$.
In~\Prop{propB+} we prove that the vector
$$v:=\gamma\bigl((T_{m}\otimes\psi)\id'\bigr)\in M(\lambda+x\xi)$$
satisfies $\HC_+(v)=ht^{-m}$ and $\cN^+ v=0$.

\subsection{}\label{v+}
Fix $\lambda\in\Lambda$ and
retain notation of~\ref{uM}.
Fix $w\in M(\lambda+x\xi)_{\lambda+x\xi-s\delta}$ 
satisfying $\cN^+  w=0$ and $z=a(-m)\in\cH_-$ ($m>0$).
In this subsection we will construct
vectors $v_{\pm}\in \breve{M}(\lambda+x\xi)^{\cN^{\pm}}$
satisfying $\HC_{\pm}(v_{\pm})=z\HC_+(w)$.

\subsubsection{}\label{Tk} 
Let $p$ be the reminder of $-m$ modulo $r$ ($0\leq p\leq r$).

Let $\fg\not=\fgl(n|n)$. Set  $T_{m+p}:u(j)\mapsto u(j-p-m)$
for $u\in\fg, j\in\mathbb{Z}$.

Let $\fg=\fgl(n|n)$. Define $T_{m+p}$
by the same formula for $j\not=0, p+m$; 
set $T_{m+p}(u(0))=\ol{u}(-p-m)$,
where $\ol{u}$ is the image of 
$u\in\fgl(n|n)^{\epsilon}$ in $\fpgl(n|n)^{\epsilon}$.
Finally, for $j=m+p$ choose any linear embedding
$\iota:\fpgl(n|n)\to \fgl(n|n)$ such that
the restriction of $\iota$ to Borel subalgebras containing
$\fn$ and $\fn_-$ are homomorphism;
set $T_{m+p}(u(m+p))=\iota(u)(0)$.

From the realization of $\fhg$ given in~\ref{XNr} one sees
that $T_{m+p}$ is well-defined linear endomorphism
of the space generated by the elements $u(j)$; this space
is the orthogonal compliment of $(\mathbb{C}K+\mathbb{C}D)$ 
in $\fhg$ with respect to the form $(-|-)$.
Observe that 
\begin{equation}\label{Tmp}
T_{m+p}([x,y])=[x,T_{m+p}(y)]=[T_{m+p}(x),y]\ \ \ \text{ if }
\wt x+\wt y\not\in\{0, (m+p)\delta\}.
\end{equation}

Fix $a(p)\in\tilde{\cH}$ such that 
$$T_{m+p}(a(p))=a(-m);$$
this element is unique unless $\fg=\fgl(n|n), p=0$.

\subsubsection{}
Set $\cH(z):=\{u\in\cH|\ [u,a(p)]=0\}$
and fix a subalgebra $\fq\subset\fhg$
such that
$$\cN^+\subset\fq\subset \cN^++\cH(z).$$
For $\fhg\not=A(2k,2l)^{(4)}$ one has $\cH(z)=\cH$ for all $z$.

\subsubsection{}
Set $N:=\mathbb{C}a(p)+\sum_{s\geq p} \cN^+_s$
and write $N=V\oplus \bigl(\sum_{s\geq m+p} \cN^+_{s}\bigr)$, where 
$$V:=\mathbb{C}a(p)+\sum_{p\leq s <m+p} 
\cN^+_s.$$
For the non-twisted case ($r=1$) one has $p=0$ and 
$N=\mathbb{C}a(0)+\sum_{s\geq 1}\fn t^s$.

Set
$$V^*:=\{f\in \Hom(N, \mathbb{C})|\ f(\sum_{s\geq m+p}
\cN^+_{s})=0\}.
$$

Both $N,\ \sum_{s\geq m+p} \cN^+_{s}$ are $(\fhh+\fq)$-modules
with respect to the adjoint action;
view  $N^*, V^*$ as $(\fhh+\fq)$-modules via the antiautomorphism $-\id$.
Observe that $V$ and $V^*$ are dual
$\fhh$-modules.
Let $a(p)^*\in V^*$ be the ``dual to $a(p)$'', that is
$$a(p)^*(a(p))=1,\ \ \ \
\forall u\in\sum_{s\geq p} {\cN}^+_s\ \ a(p)^*(u)=0.$$

\subsubsection{}\label{psiw}
Let $w\in \breve{M}(\lambda+x\xi)$ be a weight vector satisfying
$\fq_A w=0$.
By~\Lem{lemVhm} (ii) there exists a unique weight $\fq_A$-homomorphism 
$\psi:V^*_A\to M(\lambda+x\xi)$
such that $\psi(a(p)^*)=w$.

\subsubsection{}
\begin{prop}{propB+}
Let  $\gamma: \fhg_A\otimes M(\lambda+x\xi)\to M(\lambda+x\xi)$ 
be the natural map $\gamma(u\otimes w):=uw$, and  let 
$\id'\in V_A\otimes V^*_A$ 
correspond to the identity map $V_A\to V_A$.
Then 
$$v:=\gamma\bigl((T_{m+p}\otimes\psi)\id'\bigr)\in 
M(\lambda+x\xi)$$
is a weight vector such that
$\HC_+(v)=a(-m)\HC_+(w)$
and $uv=0$ for any weight element
$u\in \fq_A$ if its weight is not equal to $m\delta$.
\end{prop}
\begin{proof}
Let $B$ be a weight basis of $V\cap {\cN}^+$; then 
$B':=\{a(p)\}\cup B$
 is a weight basis of $V$. For $b\in B$ denote by 
$b^*$ the element of the dual basis $\{a(p)^*\}\cup B^*$ of $V^*$.
Then
$$v=\sum_{b\in B'} T_{m+p}(b)\psi(b^*)=
a(-m) w+\sum_{b\in B} T_{m+p}(b)\psi(b^*).$$
Notice that $T_{m+p}(b)\in {\cN}^+_-$ and so $\HC_+(v)=a(-m)\HC_+(w)$
as required.

It remains to show that $uv=0$ if 
$\wt u\not=m\delta$. Fix such element $u$.
Since $\psi$ is $\fq_A$-invariant we have

\begin{equation}
\label{himi}
u v=\sum_{b\in B'} u T_{m+p}(b)\psi(b^*)=
\sum_{b\in B'} [u,T_{m+p}(b)]\psi(b^*)+(-1)^{p(u)p(b)}T_{m+p}(b)\psi(ub^*).
\end{equation}

Since  $\wt u\not=m\delta$ the formula~(\ref{Tmp}) implies
$T_{m+p}([u,x])=[u,T_{m+p}(x)]$. Therefore
$$u v=
\sum_{b\in B'} T_{m+p}([u, b])\psi(b^*)
+(-1)^{p(u)p(b)}T_{m+p}(b)\psi(ub^*).
$$
For $b_s\in B'$ write
$$[u,b_s]=\sum_{b_j\in B'}
 c_{sj}b_j+ w_s,\ \text{ where } c_{sj}\in\mathbb{C}, 
w_s\in \sum_{t\geq m+p} {\cN}^+_t.$$
 Then $u b^*_s=(-1)^{p(u)p(b^*_s)+1}\sum c_{js} b^*_j$
and this gives $u v=\sum_s T_{m+p}(w_s)\psi(b_s^*)$.
Observe that $T_{m+p}(w_s)\in\cN^+$; since
$\psi$ is $\cN^+$-invariant we obtain
$$u v=\psi\bigl(\sum_s T_{m+p}(w_s)b_s^*\bigr).$$

We claim that $T_{m+p}(w_s)b_s^*=0$. Indeed, 
$\wt w_s=\wt u+\wt b_s$ so the weight of
$T_{m+p}(w_s)b_s^*$ is $\wt u-(m+p)\delta$.
However $(\wt u-(m+p)\delta)\not\in \Omega(V^*_A)$ because
$\Omega(V^*_A)=-\Omega(V)\subset\Omega({\cN}^-_-)\cup \{-p\delta\}$.
Hence $T_{m+p}(w_s)b_s^*=0$. This completes the proof.
\end{proof}

\subsubsection{}
\begin{rem}{v-}
Fix a subalgebra $\fq^-$ of $\fhg$ such that
$$\cN^-\subset \fq^-\subset\cN^-+\cH(z)$$
and fix a weight vector $w_-\in \breve{M}(\lambda+x\xi)^{\fq^-}$.

Substituting $\fq$ by $\fq^-$, $w$ by $w_-$, 
$N$ by  $N^-:=\mathbb{C}a(p)+\sum_{s>p} \cN^-_s$
and $V$ by $V^-:=\mathbb{C}a(p)+\sum_{p< s \leq m+p} 
\cN^+_s$,  we obtain a weight vector $v\in M(\lambda+x\xi)$
which satisfies 
$\HC_-(v)=z\HC_-(w_-)$
and $uv=0$ for any $u\in \fq^-$ provided that $\wt u\not=m\delta$.
\end{rem}

\subsection{}\label{wvlambda}
Let $v_{\lambda+x\xi}$ be the canonical generator of $M(\lambda+x\xi)$.
Retain notation of~\Prop{propB+} and construct $v$ for
$w:=v_{\lambda+x\xi}$.
The following lemmas will be used later.

\subsubsection{}
\begin{lem}{lema'}
\begin{enumerate}
\item
Take $0<j<m$. For any $d(j)\in\cH_{j}$ one has
$$\HC_+(d(j)v)=[d(j),a(-m)]v_{\lambda+x\xi};$$
\item
For any $d(m)\in\cH_{m}$ satisfying $[d(m),a(p)]=0$ one has
$$\begin{array}{l}
d(m)v=([d(m),a(-m)]-c)v_{\lambda+x\xi},\ \text{ where }\\
c:=\sum_{p\leq s<m+p} \frac{1}{2}
\str_{\fhg_s}\bigl(\ad d(-p)\circ \ad a(p)\bigr)
-\str_{\tilde{\cH}_s}\bigl(\ad d(-p)\circ \ad a(p)\bigr).
\end{array}$$
\end{enumerate}
\end{lem}
\begin{proof}
 Recall that
$$v=a(-m) v_{\lambda+x\xi}+\sum_{b\in B} T_{m+p}(b)\psi(b^*),$$
where $B\cup\{a(p)\}$ 
is a weight basis of $V=\mathbb{C}a(p)+\sum_{p\leq s<m+p}\cN^+_s$.
For (i) write
$$d(j)v=[d(j),a(-m)]v_{\lambda+x\xi}+
\sum_{b\in B} [d(j),T_{m+p}(b)]\psi(b^*)+\sum_{b\in B}
\pm T_{m+p}(b)\bigl(d(j)\psi(b^*)\bigr).$$
Since $T_{m+p}(b)\in \cN^+_-$ the term
$T_{m+p}(b)\bigl(d(j)\psi(b^*)\bigr)$ lies in the kernel of $\HC_+$.
Moreover, if $b$ is such that $[d(j),T_{m+p}(b)]\in\cN^+_-$
then $[d(j),T_{m+p}(b)]\psi(b^*)\in\Ker\HC_+$.
Otherwise $[d(j),T_{m+p}(b)]\in\cN^+$ so
$$[d(j),T_{m+p}(b)]\psi(b^*)=\psi\bigl([d(j),T_{m+p}(b)]b^*\bigr).$$
The term $[d(j),T_{m+p}(b)]b^*$ has weight $-(m+p-j)\delta$
so it is equal to zero 
(because $\Omega(V^*)\cap\mathbb{Z}\delta=\{-p\delta\}$). This proves (i).

For (ii) notice that  $d(m)V^*=0$
since there are no $\mu,\mu'\in\Omega(V^*)$ satisfying $\mu-\mu'=m\delta$.
Since $d(m)\in\fq$ this gives
$d(m)\psi(b^*)=0$ for any $b\in B$. Therefore
$$d(m)v=[d(m),a(-m)]v_{\lambda+x\xi}+
\sum_{b\in B} [d(m),T_{m+p}(b)]\psi(b^*).$$
Recall that $b\in\cN^+$ and so 
$$[d(m),T_{m+p}(b)]=[T_{m+p}(d(m)),b]=[d(-p),b]\in\cN^+.$$
Now we can rewrite the expression for $d(m)v$ in the form
$$d(m)v=[d(m),a(-m)]v_{\lambda+x\xi}+
 \psi\bigl(\sum_{b\in B} [d(-p),b]b^*\bigr).$$

Since $M(\lambda+x\xi)$ is one-dimensional,
$d(m)v=0$ if $d(m)$ and $a(-m)$ have different parities.
Note that the parities of $a(p)$ and $a(-m)$ are equal.
Therefore (ii) holds if the parities of $d(m)$ and $a(p)$ 
are not equal.
Hence we may (and will)  assume that $d(m)$
and $a(p)$ have the same parity.

The term $[d(-p),b]b^*$ has 
weight  $-p\delta$ and so $[d(-p),b]b^*=c_ba(p)^*$ for some scalar $c_b$.
Since $d(m),a(p)$ have the same parity we obtain
$$c_b=\bigl([d(-p),b]b^*\bigr)(a(p))=
(-1)^{1+p(b)} b^*\bigl([d(-p),[a(p),b]]\bigr).$$
Since $\psi(a(p)^*)=v_{\lambda+x\xi}$ one has
$$
d(m)v=Xv_{\lambda+x\xi},\ \text{where } X:=[d(m),a(-m)]+\sum_{b\in B}
(-1)^{1+p(b)} b^*\bigl([d(-p),[a(p),b]]\bigr).
$$
Notice that $B$ spans $V\cap {\cN}^+$ so
$$X=[d(m),a(-m)]
-\str_{V\cap {\cN}^+}\bigl(\ad d(-p)\circ \ad a(p)\bigr).$$

Clearly, $V\cap {\cN}^+=\sum_{p\leq s<m+p} \cN^+_s$.
Recall (see~\ref{symN}) that $\fhg$ admits 
an automorphism $\varepsilon$ satisfying 
$\varepsilon(\tilde{\cN}^+_j)=\tilde{\cN}^-_j$
which stabilizes the product $d(-p)a(p)$.
Therefore $\str_{\cN_j^+} (\ad d(-p)\circ \ad a(p))=
\str_{\cN_j^-} (\ad d(-p)\circ \ad a(p))$. Hence
$$\str_{V\cap {N}^+}\bigl(\ad d(-p)\circ \ad a(p)\bigr)=
\frac{1}{2}\str_{W}\bigl(\ad d(-p)\circ \ad a(p)\bigr)-
\str_{W\cap\tilde{\cH}}\bigl(\ad d(-p)\circ \ad a(p)\bigr),$$
where $W:=\sum_{p\leq j<m+p} \fhg_j$.
This proves (ii).
\end{proof}

\subsubsection{}
A lemma similar to~\Lem{lema'} holds for $\HC_-$.

\subsection{Proof of surjectivity in~\Prop{propz1}}
\label{HC-1z}
Take any vector weight vector 
$w_1\in \breve{M}(\lambda+x\xi)^{\cN^+}$.
Applying~\Prop{propB+} for $\fq=\cN^+$ and $w:=w_1$
we obtain 
$w_2\in \breve{M}(\lambda+x\xi)^{\cN^+}$ satisfying 
$\HC_+(w_2)=a(-m)\HC_+(w_1)$. Since the elements
of the form $a(-m)$ generate $\cS_A$, the surjectivity follows.

\subsection{Proof of injectivity in~\Prop{propz1}}\label{injHC}
Let us show that the restrictions of 
$\HC_{\pm}$ to $\breve{M}(\lambda)^{\cN^{\pm}}$ are injective.

Observe that for any $\lambda'\in\fhh$
the space $\breve{M}(\lambda')^{\cN^+}$ is orthogonal
to $\Ker\HC_-\cap\breve{M}(\lambda')$ with respect to the
Shapovalov form, see~(\ref{orth}). 
If $\lambda'$ is such that
$M(\lambda')$ is irreducible then
its  Shapovalov form is non-degenerate and thus
$$\ch \breve{M}(\lambda')^{\cN^+}\leq e^{\lambda'}\ch\cS.$$
Now take $\lambda'\in\Lambda$.
Then $\HC_+: \breve{M}(\lambda'+x\xi)^{\cN^+}\to\cS_A$ is surjective so
$\HC_+: \breve{M}(\lambda')^{\cN^+}\to\cS$
is also surjective. The above inequality gives
\begin{equation}\label{miro}
\HC_+: \breve{M}(\lambda')^{\cN^+}\iso\cS\ \ \ \text{ if 
$M(\lambda')$ is irreducible.}
\end{equation}

Set $R:=\mathbb{C}[x]$ and let $N(\lambda+x\xi)$
be a Verma $\fhg_R$-module of highest weight $\lambda+x\xi$.
Clearly,  $M(\lambda+x\xi)=N(\lambda+x\xi)\otimes_R A$.
Identify $N(\lambda+x\xi)$ with the corresponding $R$-submodule
of $M(\lambda+x\xi)$. 
Assume that $\HC_+:  \breve{M}(\lambda+x\xi)^{\cN^+}\to\cS_A$
is not injective. Then there exists a non-zero
$v\in \breve{N}(\lambda+x\xi)^{\cN^+}$ such that $\HC_+(v)=0$.
Note that the evaluation of $N(\lambda+x\xi)$ at $x=t\in\mathbb{C}$
is a Verma $\fhg$-module  $M(\lambda+t\xi)$.
The choice of $\xi$ ensures that $M(\lambda+t\xi)$ is irreducible
for $t\in\mathbb{C}\setminus Y$, where $Y$ is at most countable.
The image of any non-zero vector $v'\in N(\lambda+x\xi)$ 
 in the evaluated module
$M(\lambda+t\xi)$ is zero only for finitely many values of
$t\in\mathbb{C}$.
Hence for some $t\in \mathbb{C}\setminus Y$
the image $\ol{v}$ of $v$ in $M(\lambda+t\xi)$ is non-zero.
This contradicts to~(\ref{miro}),
since $M(\lambda+t\xi)$ is irreducible, $\cN^+\ol{v}=0$ and
 $\HC_+(\ol{v})=0$.
This establishes the injectivity and
completes the proof of~\Prop{propz1}.\qed

\section{Proof of~\ref{reform}, (A)}
\label{sectexp}
{\em In this section $\fhg$ is not of the type $A(2k,2l)^{(4)}$.}

The assertion (A) of~\ref{reform} follows 
from~\Prop{propHC} and~\Prop{propAA} below.

\subsection{Notation}\label{HC+-1}
By~\Prop{propz1} the map 
$\HC_{+}: \breve{M}(\lambda+x\xi)^{\cN^{+}}\to \cS_A$
is bijective. Let $\HC_+^{-1,\lambda}:\cS\to M(\lambda)$ be
 the evaluation
 of the inverse map $\cS_A\to \breve{M}(\lambda+x\xi)^{\cN^{+}}$
at the point $x=0$.
Clearly, $\HC_+\circ \HC_+^{-1,\lambda}=\id_{\cS}$.
Define similarly the map $\HC_-^{-1,\lambda}:\cS\to M(\lambda)$.

\subsection{}
We start  from the following lemma.

\begin{lem}{lemququ}
Take $\lambda\in\Lambda_{crit}$.
For any $u\in\cH_-$ the vectors $\HC_+^{-1,\lambda}(u),\ 
\HC_-^{-1,\lambda}(u)\in M(\lambda)$
are singular.
\end{lem}
\begin{proof}
We will prove that $\HC_+^{-1,\lambda}(u)$ is singular; the proof for
$\HC_-^{-1,\lambda}(u)$ is similar. 
We may (and will) assume that $u$ is of the form
$u=a(-m)$ for some $m>0$.

Fix any $\lambda'\in\Lambda$ and let
 $v_{\lambda'+x\xi}$ be the canonical generator of $M(\lambda'+x\xi)$.
Let $v'\in \breve{M}(\lambda'+x\xi)$ be a
vector satisfying $\cN^+v'=0,\ \ \HC_+(v')=a(-m)$.
By~\ref{propz1}, $v'$ is unique and  so it coincides with the 
vector constructing in~\Prop{propB+}
for $w=v_{\lambda'+x\xi}$. Recall that $\cH$ is commutative
(see~\ref{trfhg}) and so in~\ref{propB+} we may choose 
$\fq=\cN^++\cH$. Then~\ref{propB+} gives
$\cH_s v'=0$ for $s<m$. Clearly, $\cH_s v'=0$ for $s>m$.

Let $\ol{v}'$ be the image of $v'$  at $M(\lambda')$. By definition,
$\HC_+^{-1,\lambda'}(u)=\ol{v}'$. Notice that
$\ol{v}'\in \breve{M}(\lambda')$ and
$\cN^+\ol{v}'=\cH_s v'=0$ for $s\not=m$.
In the light of~\Lem{primcrit}, $\ol{v}'$ is singular if
$\cH_m\ol{v}=0$. Take any $d(m)\in\cH_m$.

Recall that $\cH$ is commutative. Combining~\Lem{pairing} and
\Lem{lema'} (ii) we get $d(m)v'=Xv_{\lambda'+x\xi}$, where
$$
X=m\bigl(d(m)|a(-m)\bigr)K-\frac{1}{2}\str_{W} (\ad d(-p)\circ \ad a(p)),
\ \ \ \text{ for }\ W:=\sum_{p\leq j<m+p} \fhg_j.
$$

The restriction of $(-|-)$ to $\fhg_{m}\otimes\fhg_{-m}$ is 
a non-degenerate pairing, which is invariant with respect to the adjoint 
action of $\fhg_0$. Note that the pairing
$B':\fhg_{-p}\otimes\fhg_{p}\to\mathbb{C}$
given by $B'(x|y):=\str_W (\ad x\circ \ad y)$
is also invariant with respect to the adjoint action
of $\fhg_0$. 
Identify $\fhg_{-m}$ (resp., $\fhg_{m}$)
with its image $T_{-m-p}(\fhg_{-m})\subset\fhg_{p}$ 
(resp., $T_{m+p}(\fhg_{m})\subset\fhg_{-p}$) 
via $T_{-m-p}$ (resp.,  via $T_{m+p}$). Observe that 
$\fhg_{m}$ is a simple $\fhg_0$-module (see~\cite{Kacbook}
and~\cite{vdL}, 6.10). Thus $B'$ is proportional to
the restriction of $(-|-)$; let $c$ be the coefficient
of proportionality. Note that $c$ does not depend on $\lambda'$.

We obtain $X=(d(m)|a(-m))(mK+c)$, that is
$$d(m)\ol{v}'=(d(m)|a(-m))(mK+c)v_{\lambda'}=
(d(m)|a(-m))\bigl(m(\lambda',\delta)+c\bigr)v_{\lambda'}.$$
Therefore $\ol{v}'$ is singular if $\lambda$ has level $-c/m$.
Since $\HC_+(\ol{v}')=a(-m)$, the vector $\ol{v}'$ is non-zero for all
$\lambda'\in\Lambda$.
Recall that $M(\lambda')$ is simple if $\lambda'$ has a non-critical level.
Thus  $\ol{v}'$ is not singular if the level is non-critical.
 Hence $\ol{v}'$ is singular if $\lambda'$
has the critical level.
\end{proof}

\subsection{}
\begin{prop}{propHC}
For $\lambda\in\Lambda_{crit}$
the images of $\HC_{\pm}^{-1,\lambda}:\ \cS\to M(\lambda)$
lie in $M(\lambda)^{\fhn}$.
\end{prop}
\begin{proof}
We will prove that $\HC_+^{-1,\lambda}(z)$ is singular for $z\in\cS$; 
the proof for $\HC_-^{-1,\lambda}(z)$ is similar. 

For $v\in M(\lambda+x\xi)$ denote by $\ol{v}$ its image in $M(\lambda)$.
By~\Prop{propz1} the map 
$\HC_{+}: \breve{M}(\lambda+x\xi)^{\cN^{+}}\to \cS_A$
is bijective; denote by $\HC_{+}^{-1,\lambda+x\xi}$ 
the inverse map.

For $z\in\cH_-$ the assertion is proven in~\Lem{lemququ}.
Since $\cS$ is generated by $\cH_-$, it is enough to show
that 
$$\HC_+^{-1,\lambda}(z) \text{ is singular }\ \Longrightarrow\ \ 
\HC_+^{-1,\lambda}(a(-m)z) \text{ is singular }.$$
Set $z_1:=a(-m),\ z_2:=z,\ z_3=a(-m)z$
and let $w_i:=\HC_+^{-1,\lambda+x\xi}(z_i)$ for $i=1,2,3$.
The above implication takes form:
$\fhn\ol{w}_2=0 \Longrightarrow\
\fhn\ol{w}_3=0$.
Assume that $\fhn\ol{w}_2=0$.
Let $\phi\in\End_{[\fhg,\fhg]}(M(\lambda))$ be
such that $\phi(v_{\lambda})=\ol{w}_2$ (see~\ref{lemrst}).
We claim that
\begin{equation}\label{psiw3}
\ol{w}_3=\phi(\ol{w}_1).
\end{equation}
By~\Lem{lemququ} $\ol{w}_1$ is singular so
the formula (\ref{psiw3}) implies the required implication.

Let us prove~(\ref{psiw3}). Set $\fq:=\cN^+$. 
Apply the construction described
in~\ref{v+} for $w:=v_{\lambda+x\xi}$ (resp., for $w:=w_2$)
and
denote by $\psi_1$ (resp., by $\psi_2$) the weight $\fq_A$-homomorphism
introduced in~\ref{psiw}. 
By constructions of $w_1$ (\Prop{propB+}) and of $w_3$ (\ref{HC-1z}) one has
\begin{equation}\label{psiw1}
w_1=\sum_{b\in B'} T_{m+p}(b)\psi_1(b^*),\ \ 
w_3=\sum_{b\in B'} T_{m+p}(b)\psi_2(b^*)
\end{equation}

 For $i=1,2$
denote by $\psi_i(0): V^*\to M(\lambda)$
the evaluation of $\psi_i$ at $x=0$; clearly,
$\psi_i(0)$ is a $\fq$-homomorphism.
One has
$$\psi_1(a(p)^*)=v_{\lambda+x\xi},\ \ \phi(v_{\lambda})=\ol{w}_2,\ \
\psi_2(a(p)^*)=w_2$$
so $\phi(\psi_1(0)(a(p)^*))=\psi_2(0)(a(p)^*)$.
Since both $\phi\circ\psi_1(0)$ and $\psi_2(0)$ are 
weight $\fq$-homomorphisms $V^*\to M(\lambda)$,
the uniqueness proven in~\Lem{lemVhm} (i) gives
$$\phi\circ\psi_1(0)=\psi_2(0).$$
Now~(\ref{psiw1})  gives
$$\ol{w}_3=\sum_{b\in B'} T_{m+p}(b)\phi(\psi_1(0)(b^*))=
\phi\bigr(\sum_{b\in B'} T_{m+p}(b)\psi_1(0)(b^*)\bigl)=\phi(\ol{w}_1).
$$
This establishes~(\ref{psiw3}) and completes the proof.
\end{proof}

\subsection{}
\begin{prop}{propAA}
Take $\lambda\in\Lambda_{crit}$.
For $v\in \breve{M}(\lambda+x\xi)^{\cN^{\pm}}$ one has
$v\in \cF^k(M(\lambda+x\xi))$, where $k$ is the degree
of $\HC_{\pm}(v)$.
\end{prop}
\begin{proof}
The assertions for $\cN^+$ and for $\cN^-$ have similar proofs. 
Recall that $S_A(\lambda+x\xi)$ is the Shapovalov form on
$M(\lambda+x\xi)$; denote this form by $(-,-)$.
We will prove that 

{\em for $v_+\in \breve{M}(\lambda+x\xi)^{\cN^{+}}$ and 
$v\in M(\lambda+x\xi)$
one has $(v_+,v)\in (x)^k$, where
$k$ is the degree of $\HC_+(v_+)$.}

The proof is by induction on $k$.
For $k=0$ the assertion is trivial.

Let $v_-\in M(\lambda+x\xi)^{\cN^-}$ be such that $\HC_-(v_-)=\HC_-(v)$.
Similarly to~(\ref{orth}) one has
$(v_+,\Ker\HC_-)=(\Ker\HC_+,v_-)=0$. Therefore
$$(v_+,v)=(v_+,v_-)=(\HC_+(v_+)v_{\lambda+x\xi},v_-).$$

We may (and will) assume that $z:=\HC_+(v_+)$ is of the form
$z=uw$ for $u\in\cH_-, w\in\cS$. Then
$$(\HC_+(v_+)v_{\lambda+x\xi},v_-)
=(zv_{\lambda+x\xi}, v_-)=(wv_{\lambda+x\xi},
\sigma(u)v_-).$$

Since $\sigma(u)\in\fhn$, \Prop{propHC} gives $\sigma(u)v_-=xv'$
for some $v'\in M(\lambda+x\xi)$.  Moreover, since $\sigma(u)\in\cH$
one has $\cN^-(\sigma(u)v_-)\subset \cN^-v_-=0$ so $\cN^-v'=0$. Thus
$$(v_+,v)=x (wv_{\lambda+x\xi},v'),$$
for  $v'\in \breve{M}(\lambda+x\xi)^{\cN^{-}}$. Taking
$w_+\in \breve{M}(\lambda+x\xi)^{\cN^{+}}$ such that $\HC_+(w_+)=w$
we obtain
$$(v_+,v)=x (w_+,v').$$
Notice that $w=\HC_+(w_+)$ is a monomial of degree $k-1$.
By induction hypothesis, $(w_+,v')\in (x)^{k-1}$. This completes the proof.
\end{proof}

\subsection{}\label{pfA}
Combining Propositions~\ref{propHC} and~\ref{propAA} we obtain that
$\HC_{\pm}^{-1,\lambda}(\cS^{\geq k})$ lies in $M(\lambda)^k\cap
M(\lambda)^{\fhn}$. Hence 
$\HC_{\pm}\bigl(M(\lambda)^k\cap M(\lambda)^{\fhn}\bigr)$ 
contains $\cS^{\geq k}$ and this establishes~\ref{reform} (A).

\section{Proof of~\ref{reform}, (B)}\label{sectB}
{\em  In this section $\fhg$ is not of the type $A(2k,2l)^{(4)}$
and $\lambda\in\Lambda_{crit}$.}

In this section we prove the assertion (B) of~\ref{reform}.
As it is shown in Sect.~\ref{sectorm} this completes the proof of
Theorems~\ref{thm01}, \ref{thm02}.

\subsection{Notation}
Set 
$$Z:=\frac{\ch\cU(\fhn_-)}{\ch\cS}=
\ch \cU(\cN^-_-)\ch \cU(\cN^+_-).$$

In~\ref{xy} we introduced a ring $C$
which contains, in particularly, $\ch M(\lambda),\ \ch L(\lambda)$
and $Z$. Introduce the following projections $P_{\nu}, P_{\geq\nu}: C\to C$:
$$P_{\nu}(\sum c_{\mu}e^{\mu}):=c_{\nu}e^{\nu},\ \ 
P_{>\nu}(\sum c_{\mu}e^{\mu}):=\sum_{\mu>\nu}c_{\mu}e^{\mu}.$$
Set $\ch_{>\nu} N:=P_{>\nu}(\ch N)$.

Denote by $a_{k,s}$ the multiplicity of $L(\lambda-s\delta)$
in $M(\lambda)^k$:
$$a_{k,s}:=[M(\lambda)^k:L(\lambda-s\delta)].$$

\subsection{}
Recall that $\{\cS^{\geq k}\}$ is the adic filtration on $\cS$.
We will prove by induction on $m\geq 0$ that
\begin{equation}\label{indch}\begin{array}{ll}
(i) & a_{k,j}=\dim\cS_{-j\delta}^{\geq k}\text{ for } 
j=0,\ldots,m-1\text{ and all } k;\\
(ii) & \ch_{>\lambda-m\delta}  L(\lambda)=e^{\lambda}
P_{>- m\delta} (Z).
\end{array}
\end{equation}
Clearly,~\ref{reform} (B) follows from (i) for $m=1,2,\ldots$.

\subsubsection{}
For $m=1$ (i) holds because $a_{k,0}=\delta_{k,0}=\dim\cS_0^{\geq k}$.
For (ii) recall that $[M(\lambda):L(\mu)]=0$
for $\mu>\lambda-\delta$ so
$\ch_{>\lambda-\delta}  L(\lambda)=\ch_{>\lambda-\delta}  M(\lambda)$.
Therefore 
$$\ch_{>\lambda-\delta}  L(\lambda)=
e^{\lambda} P_{>- \delta} \bigl(\cS\cdot Z\bigr)=e^{\lambda}
P_{>- \delta} (Z),$$
since $P_{>- \delta}(\cS)=1$ and $P_{\nu}(Z)=0$ for
$\nu\not\leq 0$.

\subsubsection{}
Assume that~(\ref{indch}) holds for $m=s-1$. 

Let us check (i) for $m=s$. Recall that 
$M(\lambda)^k_{\mu}=0$ for $k>>0$ and so
$a_{k,j}=0$ for $k>>j$.
Notice that $a_{k,0}=0$ for $k\geq 1$; this gives
$$\sum_{k\geq 1} \ch_{>\lambda-s\delta}  M(\lambda)^k=
\sum_{k,j\geq 1} a_{k,j} \ch_{>\lambda-s\delta}  L(\lambda-j\delta)=
\sum_{k,j\geq 1} a_{k,j}e^{-j\delta}
\ch_{>\lambda-(s-j)\delta}  L(\lambda).$$
Using the induction hypothesis we get
$$e^{-\lambda}\sum_{k\geq 1} \ch_{>\lambda-s\delta}  M(\lambda)^k=
\sum_{k,j\geq 1} a_{k,j}
e^{-j\delta}P_{>- (s-j)\delta} (Z).$$
On the other hand, the formula~(\ref{sumMS}) gives
$$e^{-\lambda}\sum_{k\geq 1} \ch_{>\lambda-s\delta}  M(\lambda)^k=
P_{>-s\delta}\bigl(Z\sum_{k\geq 1}\ch \cS^{\geq k}\bigr)=
\sum_{k,j\geq 1}\dim\cS^{\geq k}_{-j\delta}e^{-j\delta} 
P_{>- (s-j)\delta} (Z).$$

Hence
\begin{equation}\label{trtr}
\sum_{k,j\geq 1} b_{k,j}e^{-j\delta}P_{>- (s-j)\delta} (Z)=0,\  \text{ where } 
\ b_{k,j}:=a_{k,j}-\dim\cS^{\geq k}_{-j\delta}.
\end{equation}
By~\ref{reform} (A) $b_{k,j}\geq 0$ for all $k,j$.

Note that  $a_{k,j}= \cS^{\geq k}_{-j\delta}=0$ for $k>>j$.
Observe that
the term $P_{>- (s-j)\delta} (Z)=0$ iff $j\geq s$;
for $j=1,\ldots,s-1$ this term is a finite sum of the form
$\sum_i c_i e^{\mu_i}$, where $c_i\geq 0$.
Hence the sum in the left-hand side of~(\ref{trtr}) is finite.
Now combining the above observation and
the  inequalities $b_{k,j}\geq 0$,
we get $b_{k,j}=0$ for $j=1,\ldots,s-1$.
Hence $a_{k,j}=\dim\cS^{\geq k}_{-j\delta}$ for $j=1,\ldots,s-1$.
Since $a_{k,0}=\delta_{k,0}=\dim\cS^{\geq k}_{-j\delta}$
we obtain (i) for $m=s$.

\subsubsection{}
For (ii) write
$$\ch_{>\lambda-s\delta}  M(\lambda)=
\ch_{>\lambda-s\delta}  L(\lambda)+\sum_{j=1}^{s-1} a_{1,j} 
\ch_{>\lambda-s\delta} L(\lambda-j\delta).$$
Using (i) we obtain
$$\begin{array}{rl}
\ch_{>\lambda-s\delta}  L(\lambda)&=
\ch_{>\lambda-s\delta}  M(\lambda)-\displaystyle\sum_{j=1}^{s-1}
\dim\cS_{-j\delta}\ch_{>\lambda-s\delta} L(\lambda-j\delta)\\
& =
e^{\lambda}P_{>-s\delta}(\ch\cS\cdot Z) 
-\displaystyle\sum_{j=1}^{s-1}\dim\cS_{-j\delta} e^{-j\delta}
\ch_{>\lambda-(s-j)\delta} L(\lambda).\end{array}$$
Finally, by the induction hypothesis,
$$e^{-\lambda}\ch_{>\lambda-s\delta}  L(\lambda)=
P_{>-s\delta}(\ch\cS\cdot Z) 
-\sum_{j=1}^{s-1}\dim\cS_{-j\delta} e^{-j\delta}
P_{>-(s-j)\delta}(Z)=P_{>-s\delta}(Z),$$
and this establishes (ii).

\section{Vanishing lemmas}\label{sectvan}
\subsection{}
The following lemma is proven in~\cite{GS}, 8.1.

\begin{lem}{chmlem1}
Let $\fm$ be a subalgebra of $\fhn$. Assume that
$\lambda,\mu\in\fhh^*$ are such that 
$$L(\lambda)_{\mu}\subset \sigma(\fq)L(\lambda),\ \ 
M'(\lambda)_{\mu}\subset \sigma(\fm)M'(\lambda)$$
then
$$H^r(\fm, L(\lambda))_{\mu}=0\ \text{ for } r=0,1.$$
\end{lem}

\subsection{}
Let $\fm$ be a Lie superalgebra which carries 
a $\mathbb{Z}^k_{\geq 0}$-grading
$\fm=\oplus_{\nu\in\mathbb{Z}^k_{\geq 0}}\fm_{\nu}$ 
with  finite-dimensional
homogeneous components.  Set 
$$R:=\mathbb{C}[x],\ \ \ A:=\mathbb{C}[x]_{(x)},$$
that is $A$ is the localization of $R$ by the ideal generated by $x$.

Set $\fm_R:=\fm\otimes_{\mathbb{C}} R$.
Let $M=\oplus_{\nu\in\mathbb{Z}^k_{\leq 0}}M_{\nu}$ 
be a graded $\fm_R$-module whose
 homogeneous components $M_{\nu}$ are free $R$-modules of finite rank.

For each $c\in\mathbb{C}$ view  $M(c):=M/(t-c)M$ as an $\fm$-module. 
Set $\fm_A:=\fm\otimes_R A$, $M_A:=M\otimes_R A$
and view $M_A$ as an $\fm_A$-module. In the following lemma
$\fm_A$ is considered as a Lie superalgebra over $A$.

\subsubsection{}
\begin{lem}{lemcohoR}
One has 
\begin{enumerate}
\item
$H^0(\fm,M(0))_{\nu}=0\ 
\Longrightarrow\ H^0(\fm_A,M_A)_{\nu}=0$.
\item Suppose that
 $H^0(\fm,M(c))_{\nu}=H^1(\fm,M(c))_{\nu}=0$ for $c\in\{0\}\cup S$,
where $S$ is an infinite subset of $\mathbb{C}$.
Then $H^1(\fm_A,M_A)_{\nu}=0$.
\end{enumerate}
\end{lem}
\begin{proof}
(i) If $H^0(\fm_A,M)_{\nu}\not=0$ then there exists $v\in M_{\nu}$
such that $\fm v=0$. Write $v=t^kw$ for
 $w\not\in (t)M$. Then  $\fm w=0$ and the image of
$w$ in $M(0)=M/(t)M$ is non-zero; this contradicts to
$H^0(\fm,M(0))_{\nu}=0$.

(ii)  Let $\phi:\fm_{A}\to M_A$ be an $A$-linear map 
satisfying the conditions $\phi(\fm_{A,\mu})=M_{\mu+\nu}$ and
$\phi([u_1,u_2])=u_1\phi(u_2)-(-1)^{p(u_1)p(u_2)}u_2\phi(u_1)$.
The required equality $H^1(\fm_A,M_A)_{\nu}=0$ means that
$M_{A,\nu}$ contains a vector $v$ such that 
$uv=\phi(u)$ for any $u\in\fm$.

Observe that $\fm_{A,\mu} M_{A,\nu}=0$
and $\phi(\fm_{A,\mu})=0$ if
$-(\mu+\nu)\not\in\mathbb{Z}_{\geq 0}^k$
and that $\sum_{\mu:-(\mu+\nu)\in\mathbb{Z}_{\geq 0}^k }\fm_{A,\mu}$
has a finite rank over $A$.
Let $k$ be the rank of $M_{\nu}$ over $R$ and
$B$ be a basis of $M_{\nu}$. Write
$v=\sum_{b\in B} x_b b$ with $x_b\in A$. The condition $uv=\phi(u)$
for all $u\in\fm$ can be written as a finite system of linear equations
on $\{x_b\}_{b\in B}$ (the finiteness follows from the above observation),
i.e. $DX=Y$, where $X=(x_b)_{b\in B}$
and $D$ is a matrix, $Y$ is a vector  with entries in $A$.
Let $S_1\subset S\cup\{0\}$ consist of $c$'s such that the entries of 
$D$ and $Y$ have
no poles at $c$. Note that $S_1$ is infinite and contains $0$.

For $c\in S_1$ let $D(c),Y(c)$ be 
the evaluations of respectively $D$ and $Y$
at $t=c$. We claim that the system $D(c)X'=Y(c)$ has a unique
solution for any $c\in S_1$, that is
\begin{equation}\label{Dc}
\forall c\in  S_1\ \  \exists !\ X'\in\mathbb{C}^k\ \text{s.t.}\ 
D(c)X'=Y(c).
\end{equation} 
By~\Lem{lemlast}, (\ref{Dc}) implies that $DX=Y$ has a unique solution
and so establishes (ii).

It remains to verify~(\ref{Dc}).
For $c\in S_1$
let ${\phi}(c): \fm\to M(c)$ be the evaluation of $\phi$ at $t=c$.
The assumption $H^0(\fm,M(c))_{\nu}=H^1(\fm,M(c))_{\nu}=0$ means that
 there exists a unique vector 
$w\in M(c)_{\nu}$ such that ${\phi}(c)(u)=uw$ for all  $u\in \fm$.
Let $B(c)$ be the image $B$ in $M(c)$.
Then $w=\sum_{b'\in B(c)} x'_{b'} b'$ with $x'_{b'}\in \mathbb{C}$, where
 $X'=(x'_{b'})_{b'\in B(c)}$ satisfies the system $D(c)X'=Y(c)$
and this establishes~(\ref{Dc}) and completes the proof.
\end{proof}

\subsubsection{}
Set $F:=\Fract \mathbb{C}[x]$. 
For a matrix $D$ over $F$ and $c\in\mathbb{C}$ 
we say that ``$D$ has no poles at $t=c$'' if
the entries of $D$ have no poles at $c$; in this case we denote by
$D(c)$ its evaluation at $c\in\mathbb{C}$.

\begin{lem}{lemlast}
Take $D\in\Mat_{m\times k}(F), Y\in F^m$. Assume that
there exists  an infinite set $S\subset\mathbb{C}$ such that 
$D,Y$ have no poles in $S$ and 
$$\forall c\in  S\ \  \exists !\ X'\in\mathbb{C}^k\ \text{s.t.}\ 
D(c)X'=Y(c).$$
Then there exists a unique $X\in F^k$ such that $DX=Y$; moreover, $X$
has no poles in $S$.
\end{lem}
\begin{proof}
Fix $s\in S$. Since $D(s)X'=Y(s)$ has a unique solution,
$D(s)$ has a square submatrix $C'$ 
of size $k$ such that $\det C'\not=0$. Let $C$ be 
the corresponding submatrix of $D$ and $Z$ be
the corresponding   sub-column of $Y$ 
(i.e. $Z$ is the transpose of $(y_{i_1},\ldots,y_{i_k})$, where
$C=(d_{ij})_{i_1,\ldots,i_k; j=1,\ldots,k}$).
Clearly,  $C'=C(s)$, 
and $(\det C)(s)=\det C(s)$ so
$(\det C)(s)\not=0$. In particular, $\det C\not=0$ so
$C$ is invertible and $DX=Y$ has at most one solution.

Let us show that $DX=Y$ for
$X:=C^{-1}Z$ and that $X$ has no poles in $S$.

Let $S_1:=\{c\in S|\ (\det C)(c)\not=0\}$. 
Note that $S_1$ is infinite and contains $s$;
$X$ has  no poles in $S_1$.
Take $c\in S_1$. 
Then $\det C(c)\not=0$ so $C(c)$ is a non-zero minor of the maximal size
in the system $D(c)X'=Y(c)$. Since
this system has a unique solution, this solution
is $C(c)^{-1}Z(c)=X(c)$. Hence $D(c)X(c)=Y(c)$
for any $c\in S_1$.  Since $S_1$ is infinite and the entries
of $D(c),X(c)$ and $Y(c)$ are rational functions in $c$,
we conclude that $DX=Y$. 
Hence a unique solution of the system $DX=Y$ has no poles at $s$.
Since $s$ is an arbitrary element of $S$ this solution has no poles
in $S$.
\end{proof}


\end{document}